# CONTINUUM TREE ASYMPTOTICS OF DISCRETE FRAGMENTATIONS AND APPLICATIONS TO PHYLOGENETIC MODELS

By Bénédicte Haas, Grégory Miermont, Jim Pitman
and Matthias Winkel[1]

*Université de Paris Dauphine, Université de Paris Sud, University of
California, Berkeley and University of Oxford*

Given any regularly varying dislocation measure, we identify a natural self-similar fragmentation tree as scaling limit of discrete fragmentation trees with unit edge lengths. As an application, we obtain continuum random tree limits of Aldous's beta-splitting models and Ford's alpha models for phylogenetic trees. This confirms in a strong way that the whole trees grow at the same speed as the mean height of a randomly chosen leaf.

**1. Introduction.** For a number of years, there has been an increased interest in random tree models, both in the mathematical literature and in applied sciences such as genetics. Fundamental classes of trees are trees with $n$ leaves and no degree-2 nodes. Denote by $\mathbb{T}_n^\circ$ the space of such trees, which can be made mathematically precise as a space of connected acyclic graphs with $n+1$ degree-1 vertices, one of which is distinguished as the root. Also, denote by $\mathbb{T}_n$ the space of such trees where the other $n$ degree-1 vertices (the leaves) are labeled $1, \ldots, n$. Such trees are called *cladograms* in the genetics literature, up to trivial differences and an extension. Here, the trees are planted, that is, the root has only one neighbor, and they are not necessarily binary, as is usually assumed in the phylogenetics literature. The only edge adjacent to the root is called the *root-edge*.

A class of probability distributions on $\mathbb{T}_n^\circ$ or $\mathbb{T}_n$ can be specified by a procedure called *Markov branching* [1]: $P_1^\circ$ is the unique distribution on the singleton $\mathbb{T}_1^\circ$. Recursively, $P_n^\circ$ is the distribution of a random tree $T_n^\circ$, where

Received April 2006; revised November 2007.
[1]Supported in part by EPSRC Grant GR/T26368/01 and NSF Grant DMS-0405779.

the unique branch point neighboring the root connects $r \geq 2$ subtrees with $k_1 \geq \cdots \geq k_r \geq 1$ leaves, respectively, $k_1 + \cdots + k_r = n$, with some probability $q_n(k_1, \ldots, k_r)$ (so that $q_n$ is a probability distribution on the set of partitions of the integer $n$). Given the branching into sizes $k_1, \ldots, k_r$, these subtrees are independent, with distributions $P^\circ_{k_r}$. Finally, define $P_n$ as the distribution of the random tree $T^\circ_n$, equipped with leaf labels uniform among all possible labelings with $\{1, \ldots, n\}$. Therefore, all Markov branching models on $\mathbb{T}_n$, thus defined, have exchangeable leaf labels.

Important models in phylogenetics such as the Yule, uniform and comb models, and, more generally, Aldous's beta-splitting models [1, 5] and Ford's alpha models [21] have the Markov branching property (see [1, 21] for references to the phylogenetics literature). They also have a property of sampling consistency, that is, the subtree of $T_n \sim P_n$ generated by the leaves labeled $1, \ldots, n-1$ has distribution $P_{n-1}$.

By a standard argument using Kolmogorov's extension theorem, for sampling consistent $(P_n)_{n \geq 1}$, one can consider a *strongly sampling consistent* sequence $(T^\circ_n)_{n \geq 1}$ [resp. $(T_n)_{n \geq 1}$] defined on some probability space, in the (stronger) sense that $T^\circ_{n-1}$ is the subtree of $T^\circ_n \sim P^\circ_n$ generated by $n-1$ leaves chosen uniformly at random (resp. leaves $1, \ldots, n-1$) for all $n \geq 2$.

The recursive definition of $P_n$ is due to Aldous [1] in the binary case, where $q_n$ is supported by partitions of $n$ of the form $(n-k, k), 1 \leq k \leq n/2$, for all $n$. Not all Markov branching models are sampling consistent [e.g., a splitting rule for which $q_4(2,2) = 1$ cannot be sampling consistent] and Aldous leaves as an open problem the characterization of sampling consistent Markov branching models (on cladograms). Ford [21] gives an answer in terms of a certain consistency condition that has to be satisfied by the associated (binary) splitting rules $q_n$. See also (14) for the general nonbinary case.

A more explicit answer in the form of an integral representation, also for nonbinary models, can be obtained from Bertoin's study of homogeneous fragmentation processes [8]. In the present paper, we interpret sampling consistent Markov branching models as trees associated with (discrete) fragmentations, where $T^\circ_n \sim P^\circ_n$ describes the fragmentation of an initial mass of size $n$ (or of the set $\{1, \ldots, n\}$ for $T_n \sim P_n$), first into blocks of sizes $k_1, \ldots, k_r$ and then of each block, recursively, until all blocks have unit size. These branching models can be characterized in terms of homogeneous fragmentation processes, as follows.

A homogeneous fragmentation process is a continuous-time continuous-mass analog of the above discrete fragmentations. The most intuitive class is that of mass fragmentation processes, that is, right-continuous Markov processes $(F_t)_{t \geq 0}$ in

$$\mathcal{S}^\downarrow = \left\{ (s_i)_{i \geq 1} : s_1 \geq s_2 \geq \cdots \geq 0, \sum_{i \geq 1} s_i \leq 1 \right\},$$

whose transition kernels have the property that given state $\mathbf{s} = (s_i)_{i \geq 1}$ at time $u$, each fragment of mass $s_i$ evolves independently and with distribution identical to $(s_i F_t)_{t \geq 0}$. More precisely, for each $t \geq 0$, $F_{u+t}$ can be written as the decreasing rearrangement of masses of all fragments of $s_i F_t$, $i \geq 1$. Bertoin has shown that the distribution of such a process is determined by an *erosion coefficient* $c \in \mathbb{R}_+$ and a *dislocation measure* $\nu$ on $\mathcal{S}^\downarrow$ which satisfies

$$(1) \qquad \nu(\{(1, 0, \ldots)\}) = 0 \quad \text{and} \quad \int_{\mathcal{S}^\downarrow} (1 - s_1) \nu(d\mathbf{s}) < \infty.$$

In the sequel, for $\mathbf{s} \in \mathcal{S}^\downarrow$, we let $s_0 = 1 - \sum_{i \geq 1} s_i \in [0, 1]$.

THEOREM 1. *Sampling consistent splitting rules $(q_n, n \geq 2)$ are all of the following form: if $(k_1, \ldots, k_r)$ is a partition of $n$ with $r \geq 2$ parts, of which exactly $m_i \geq 0$ parts are equal to $i$, $1 \leq i \leq n$, then*

$$(2) \quad q_n^{(c,\nu)}(k_1, \ldots, k_r)$$
$$:= \frac{C_{k_1,\ldots,k_r}}{Z_n} \left( nc \mathbf{1}_{\{r=2, k_2=1\}} + \int_{\mathcal{S}^\downarrow} \sum_{l=0}^{m_1} \binom{m_1}{l} \sum_{\substack{i_1,\ldots,i_{r-l} \geq 1 \\ \text{distinct}}} s_0^l \prod_{j=1}^{r-l} s_{i_j}^{k_j} \nu(d\mathbf{s}) \right)$$

*for some pair $(c, \nu)$, where $c \geq 0$ and $\nu$ satisfies (1). Here, $C_{k_1,\ldots,k_r}$ is a combinatorial factor and $Z_n$ the normalization constant, as follows:*

$$(3) \quad C_{k_1,\ldots,k_r} = \frac{n!}{k_1! \ldots k_r! m_1! \ldots m_n!}, \qquad Z_n = nc + \int_{\mathcal{S}^\downarrow} \left( 1 - \sum_{i \geq 1} s_i^n \right) \nu(d\mathbf{s}).$$

*Moreover, one has $q_n^{(c,\nu)} = q_n^{(c',\nu')}$ for every $n \geq 2$ if and only if $(c', \nu') = (Kc, K\nu)$ for some $K > 0$.*

The intuitive meaning of (2) is that $(s_i)_{i \geq 0}$ is chosen according to $\nu$ and an independent sample from $(s_i)_{i \geq 0}$ of size $n$ is taken, jointly conditioned not to have all sample points in one fragment $s_i$, $i \geq 1$. The term $s_0$ is special in that each of the $l$ sample points in $s_0$ is considered a singleton. We note that

$$C_{k_1,\ldots,k_r} \binom{m_1}{l} = \binom{n}{l, k_1, \ldots, k_{r-l}} \left( (m_1 - l)! \prod_{i=2}^n m_i! \right)^{-1},$$

where $\binom{n}{l, k_1, \ldots, k_{r-l}}$ is the number of permutations with the same given frequencies $l, k_1, \ldots, k_{r-l}$. The allocation of indices $i_j$ to box sizes $k_j$ is such that $(m_1 - l)! \prod_{i=2}^n m_i!$ sequences $(i_1, \ldots, i_{r-l})$ lead to the same configuration $\{(i_1, k_1), \ldots, (i_{r-l}, k_{r-l})\}$ and hence contribute to the coefficient of the same monomial $s_0^l s_{i_1}^{k_1} \cdots s_{i_{r-l}}^{k_{r-l}}$.

A homogeneous fragmentation process for which $c = 0$, is said to have *no erosion*. Also, a dislocation measure $\nu$ is said to be *conservative* if

$$\nu(\{\mathbf{s} \in \mathcal{S}^\downarrow : s_0 > 0\}) = 0. \tag{4}$$

The conditions that $c = 0$ and $\nu$ is conservative are equivalent to the *dust-free* property of the associated homogeneous fragmentation $(F_t)_{t \geq 0}$, namely that the terms of $F_t$ sum to 1 for any $t \geq 0$, a.s. Under these conditions, formula (2) takes the much simpler form

$$\begin{aligned}
&q_n^{(0,\nu)}(k_1, \ldots, k_r) \\
&= \frac{1}{Z_n} \binom{n}{k_1, \ldots, k_r} \left( \prod_{i=1}^n m_i! \right)^{-1} \int_{\mathcal{S}^\downarrow} \sum_{\substack{i_1, \ldots, i_r \geq 1 \\ \text{distinct}}} \prod_{j=1}^r s_{i_j}^{k_j} \nu(d\mathbf{s}).
\end{aligned} \tag{5}$$

In [9], Bertoin introduced self-similar fragmentation processes $(F_t^{(a)})_{t \geq 0}$ in $\mathcal{S}^\downarrow$ with parameter $a \in \mathbb{R}$: given state $\mathbf{s}$ at time $u$, the evolution of each fragment of mass $s_i$ is independent and distributed as $(s_i F_{s_i^a t}^{(a)})_{t \geq 0}$. Once $a$ is fixed, such processes are in one-to-one correspondence with homogeneous fragmentations and are hence characterized by an erosion parameter $c$ and a dislocation measure $\nu$ satisfying (1). In the sequel, we will only deal with negative index $a < 0$ and write $\gamma = -a$.

For $\gamma > 0$, $c = 0$ and conservative dislocation measures $\nu$, associated fragmentation trees $\mathcal{T}_{(\gamma,\nu)}$ have been studied in [29] using Aldous's continuum random tree (CRT) formalism of trees as subsets of $l_1 = l_1(\mathbb{N})$ (cf. [2, 3, 4]). Alternative tree representations have been developed and we shall here use abstract $\mathbb{R}$-trees as introduced for use in probability by Evans and co-authors [18, 19, 20] (see also [27]). Following these references, the space of $\mathbb{R}$-trees will be endowed with the Gromov–Hausdorff metric, which provides a notion of convergence for these abstract spaces. All the necessary concepts are discussed in Section 3.3.

Under the regular variation condition

$$\nu(s_1 \leq 1 - \varepsilon) = \varepsilon^{-\gamma_\nu} \ell\left(\frac{1}{\varepsilon}\right) \tag{6}$$

for some $\gamma_\nu \in (0,1)$ and a function $x \to \ell(x)$ slowly varying as $x \to \infty$, the case $\gamma = \gamma_\nu$ is special. Under the further regularity condition

$$\int_{\mathcal{S}^\downarrow} \sum_{i \geq 2} s_i |\ln(s_i)|^\rho \nu(d\mathbf{s}) < \infty \tag{7}$$

for some $\rho > 0$ [this is satisfied, e.g., if $\nu(s_{r+1} > 0) = 0$ for some $r > 0$], our main theorem identifies the $\gamma_\nu$-self-similar fragmentation tree as a scaling limit of discrete fragmentation trees associated with a (homogeneous) fragmentation process or, equivalently, by Theorem 1, associated with sampling consistent Markov branching models.

THEOREM 2. *Let $\nu$ be a conservative dislocation measure satisfying* (6) *and* (7) *and* $(T_n^\circ)_{n\geq 1}$ *a strongly sampling consistent family of discrete fragmentation trees* $T_n^\circ \sim P_n^\circ$ *as associated via* (5). *If we consider* $T_n^\circ$ *as a random $\mathbb{R}$-tree (with unit edge lengths), then there is the convergence in probability*

$$\frac{T_n^\circ}{n^{\gamma_\nu}\ell(n)\Gamma(1-\gamma_\nu)} \xrightarrow[n\to\infty]{(p)} \mathcal{T}_{(\gamma_\nu,\nu)}$$

*with respect to the Gromov–Hausdorff metric, where $\mathcal{T}_{(\gamma_\nu,\nu)}$ is a $\gamma_\nu$-self-similar fragmentation tree with dislocation measure $\nu$, defined as a random $\mathbb{R}$-tree on the same probability space that supports $(T_n^\circ)_{n\geq 1}$.*

Note that we obtain a convergence of objects with constant edge lengths to objects which, heuristically, may be expected to have "shorter" edge lengths close to the leaves, where the fragmentation rate goes to infinity. Here, we find that all sufficiently regular dislocation measures $\nu$ have an intrinsic self-similarity parameter $\gamma_\nu$, which gives a natural scale for the whole tree. As an application, we obtain limiting continuum random trees for alpha and beta-splitting models. In [1], Aldous introduced a wide class of sampling consistent binary Markov branching models, via splitting rules $q_n(n-k,k)$, $1 \leq k \leq \lfloor n/2 \rfloor$, $n \geq 2$, which he symmetrized to model a planar order so that $\tilde{q}_n(k) = \tilde{q}_n(n-k) = \frac{1}{2}q_n(n-k,k)$ for $1 \leq k < n/2$ and $\tilde{q}_n(k) = q_n(n-k,k)$ if $n = 2k$ is even. That is, $\tilde{q}_n$ is the distribution of a block selected by a fair coin toss from the split of a block of size $n$. He then studied in more detail the one-parameter family

$$\tilde{q}_n(k) = \frac{1}{Z_n^{(\beta)}} \int_0^1 \binom{n}{k} x^{k+\beta}(1-x)^{n-k+\beta}\,dx$$
$$= \frac{1}{Z_n^{(\beta)}} \binom{n}{k} \frac{\Gamma(\beta+k+1)\Gamma(\beta+n-k+1)}{\Gamma(n+2\beta+2)},$$

where $\beta > -2$. This beta-splitting model satisfies the conditions of Theorem 2 for $-2 < \beta < -1$ with $\gamma = -\beta - 1$. As an important case, when $\beta = -3/2$, the tree $T_n$ is uniform on the binary trees of $\mathbb{T}_n$. Thus, we reobtain Aldous's theorem [2], stating that the scaling limit of uniform random variables on $\mathbb{T}_n$ is the celebrated Brownian continuum random tree, with self-similarity index $-1/2$. This will be discussed in more detail in Sections 2.4 and 5.1.

In [21], Ford studied a model based on a simple sequential construction as follows. The tree $T_1^\circ$ is the unique single-leaf tree in $\mathbb{T}_1^\circ$. Given $T_n^\circ$, choose one of its edges according to a weight $1 - \alpha$ for an edge between a leaf and another vertex, and a weight $\alpha$ for an edge between two other vertices. Split the edge in two, introduce a new vertex between the two edges and add another edge from the new vertex to a new leaf. The new random tree is called $T_{n+1}^\circ$. It is implicit in the work of Ford that this model satisfies the conditions of Theorem 2 if $0 < \alpha < 1$, with $\gamma = \alpha$, so there is a CRT limit. We discuss this in more detail in Sections 5.2–5.3.

Section 2 carefully introduces the discrete framework and establishes the characterization of sampling consistent splitting rules in terms of dislocation measures of homogeneous fragmentation processes (Theorem 1). Section 3 introduces the fragmentation CRTs that appear as limits in Theorem 2. Section 4 establishes the proof of Theorem 2. Specifically, we check finite-dimensional convergence for Theorem 2 (Proposition 7), provide a tightness estimate (Proposition 9) that allows the extension of finite-dimensional convergence to convergence in the Gromov–Hausdorff sense (Section 4.2) and give a version of Theorem 2 for convergence of height functions (Theorem 15), allowing a planar order and a mass measure to be carried over to the limiting CRT. The latter convergence was conjectured by Aldous [1] in the special case of the beta-splitting models. Section 5 concludes with applications to alpha, beta-splitting and stable trees.

**2. Markov branching models and discrete fragmentations trees.** Our goal in this section is to identify the sampling consistent Markov branching models on labeled trees with laws of trees that are naturally associated with homogeneous fragmentations. As first discussed in Bertoin [8], a convenient way to study homogeneous fragmentation processes is to use a "discretization of space." This amounts to considering processes that take their values in the set $\mathcal{P}$ of partitions of the set $\mathbb{N} = \{1, 2, \ldots\}$, rather than in $\mathcal{S}^\downarrow$. To study these, we need some terminology and notation.

2.1. *Partitions.* For $B \subseteq \mathbb{N}$, we let $\mathcal{P}_B$ denote the set of partitions of $B$ into disjoint nonempty *blocks*, so $\mathcal{P} = \mathcal{P}_\mathbb{N}$. For $\pi \in \mathcal{P}_B$, we write $B' \in \pi$ to indicate that $B'$ is a block of $\pi$ and $i \overset{\pi}{\sim} j$ to indicate that $i, j \in B$ belong to the same block of $\pi$. We let $\pi_1, \pi_2, \ldots$ be the blocks of $\pi$ ranked by order of least element, so $\pi_1$ is the block containing the least element of $B$, $\pi_2$ is the block containing the least element not in $\pi_1$ and so on, with the convention that $\pi_k = \varnothing$ if $\pi$ has strictly fewer than $k$ blocks. Thus, any element $\pi$ of $\mathcal{P}_B$ can be represented as a sequence $(\pi_1, \pi_2, \ldots)$, which might eventually be constant, equal to $\varnothing$. We also let $\pi_{(i)}$ denote the block of $\pi$ that contains the integer $i \in B$. If $\pi \in \mathcal{P}_B$ and $B' \subseteq \mathbb{N}$, we let $\pi|_{B'} = B' \cap \pi$ be the partition of $B' \cap B$ obtained by restricting $\pi$ to the elements of $B' \cap B$. We let $\pi|_n = \pi|_{[n]}$ for every $n \geq 1$, where $[n] = \{1, \ldots, n\}$. By convention, we let $\mathbf{1}_B$ be the trivial partition $(B, \varnothing, \ldots)$ of $B$ and $\mathbf{0}_B = (\{i_1\}, \{i_2\}, \ldots)$ be the partition of $B$ into singletons, where $i_1 < i_2 < \cdots$ is the ranked list of elements of $B$.

In the sequel, the set $\mathcal{P}$ will be endowed with the distance $\Delta(\pi, \pi') = 2^{-N(\pi,\pi')}$, where $N(\pi, \pi') = \sup\{n \geq 1 : \pi|_n = \pi'|_n\} \in \mathbb{N} \cup \{0, \infty\}$, and the associated Borel $\sigma$-algebra.

We say that a partition $\pi \in \mathcal{P}_B$ is *finer than* $\pi' \in \mathcal{P}_B$, and write $\pi \preceq \pi'$, if any block of $\pi$ is included in some block of $\pi'$. This defines a partial order $\preceq$ on $\mathcal{P}_B$. A process or a sequence with values in $\mathcal{P}_B$ is called *refining* if it is decreasing for this partial order.

2.2. *Trees.* There is a natural relation between trees with labeled leaves and refining partition-valued processes. Write $B \subset_f \mathbb{N}$ if $B$ is a finite subset of $\mathbb{N}$. For $B \subset_f \mathbb{N}$ with $n$ elements, we let $\mathbb{T}_B$ denote the set of $\mathbf{t}$, where each $\mathbf{t}$ is a collection of subsets of $B$ and also contains $\text{ROOT} \in \mathbf{t}$, such that:

- $B \in \mathbf{t}$—we call $B$ the *common ancestor* in $\mathbf{t}$;
- $\{i\} \in \mathbf{t}$ for all $i \in B$—we call $\{i\}$, $i \in B$, the *leaves* of $\mathbf{t}$;
- for all $A, C \in \mathbf{t}$, either $A \cap C = \varnothing$, or $A \subseteq C$ or $C \subseteq A$.

$A \in \mathbf{t}$ is called a *descendant* of $C \in \mathbf{t}$ if $A \subset C$, and $C$ is then called an *ancestor* of $A$. A set $A$ is called a *child* of $C$ and $C$ is called the *parent* of $A$ if $A \subset C$ and for all $D \in \mathbf{t}$ with $A \subseteq D \subseteq C$ either $A = D$ or $D = C$. If we equip $\mathbf{t}$ with the parent-child relation and also relate ROOT with $B \in \mathbf{t}$, then $\mathbf{t}$ is a rooted connected acyclic graph so that $\mathbb{T}_{[n]}$ can be identified with $\mathbb{T}_n$ in the notation of the [Introduction](#).

For $\mathbf{t} \in \mathbb{T}_B$ and $C \in \mathbf{t}$ with $k$ children $A_1, \ldots, A_k \in \mathbf{t}$, $(A_1, \ldots, A_k)$ is a partition of $C$. We can define the subtrees $\mathbf{t}_{A_1}, \ldots, \mathbf{t}_{A_k}$ pending from $C$ as $\mathbf{t}_{A_i} = \{\text{ROOT}\} \cup \{A \in \mathbf{t} : A \subseteq A_i\}$. Then $\mathbf{t}_{A_i}$ is an element of $\mathbb{T}_{A_i}$ for $1 \leq i \leq k$. Conversely, for any finite sequence of trees $\mathbf{t}_1 \in \mathbb{T}_{B_1}, \ldots, \mathbf{t}_k \in \mathbb{T}_{B_k}$, where $B_1, \ldots, B_k$ are the nonempty blocks of a partition of some $B \subset_f \mathbb{N}$, we define $\langle \mathbf{t}_1, \ldots, \mathbf{t}_k \rangle = \{B\} \cup \bigcup_{i=1}^k \mathbf{t}_i \in \mathbb{T}_B$.

DEFINITION 1. Let $(\pi(t), t \geq 0)$ take values in $\mathcal{P}_B$ for some $B \subset_f \mathbb{N}$ and be refining. Assume, further, that $\pi(0) = \mathbf{1}_B$ and $\pi(t) = \mathbf{0}_B$ for some finite $t > 0$. We define the associated fragmentation tree to be $\mathbf{t}_\pi = \{\text{ROOT}\} \cup \{A \subseteq B : A \in \pi(t) \text{ for some } t \geq 0\}$.

A similar association can be made for refining sequences $(\pi(0), \pi(1), \ldots, \pi(m))$ of partitions of some $B \subset_f \mathbb{N}$ starting at $\pi(0) = \mathbf{1}_B$ and ending at $\pi(m) = \mathbf{0}_B$.

2.3. *Homogeneous fragmentations.* If $\Pi$ is a random variable with values in $\mathcal{P}_B$, then we say that $\Pi$ is *exchangeable* if its law is invariant under the natural action of the permutations of $B$. Similarly, a $\mathcal{P}_B$-valued process $(\Pi(t), t \geq 0)$ is exchangeable if its law is invariant under the action of permutations of $B$.

DEFINITION 2. Let $B \subset \mathbb{N}$ and consider a $\mathcal{P}_B$-valued Markov process $(\Pi(t), t \geq 0)$. We assume that for every $t, t' \geq 0$, the distribution of $\Pi(t + t')$, given $\Pi(t) = \pi$, is the same as that of the random partition whose blocks are given by

$$\Pi_i(t) \cap \pi_j^i, \qquad i, j \geq 1,$$

where $\pi^1, \pi^2, \ldots$ is an i.i.d. sequence of exchangeable partitions of $\mathbb{N}$. Then the process $\Pi$ is called a homogeneous fragmentation of $B$.

When a homogeneous fragmentation of $B$ starts from the trivial partition $\mathbf{1}_B$ of $B$, we say that the process is *standard*. We will also assume nondegeneracy of the process, namely that it is not constant a.s. It is then elementary from the definition that (nondegenerate) homogeneous fragmentation processes are refining processes whose blocks all decrease to singletons. In view of the preceding section (Definition 1), this allows us to introduce the following definition.

DEFINITION 3. Let $(\Pi(t), t \geq 0)$ be a standard homogeneous fragmentation of some *finite* $B \subset \mathbb{N}$. The tree $T_B := \mathbf{t}_\Pi \in \mathbb{T}_B$ is called the *discrete fragmentation tree* associated with $\Pi$.

As argued by Bertoin, a $\mathcal{P}$-valued process $\Pi$ is a homogeneous fragmentation if and only if its restrictions to $[n]$ are homogeneous fragmentations of $[n]$, $n \geq 1$. In other words, homogeneous fragmentations of $\mathbb{N}$ are the same as consistent families of homogeneous fragmentations of $[n]$, $n \geq 1$. Obviously, this amounts to a consistency property for the associated sequence $T_{[n]}$, $n \geq 1$, of fragmentation trees, namely that $T_{[n]}$ is the tree obtained from $T_{[n+1]}$ by removing the leaf with label $n + 1$ (and the internal vertex if it has only two other neighbors, that will then be connected by a direct edge instead). We claim that the laws $(P_n, n \geq 1)$ associated with sampling consistent splitting rules as explained in the Introduction are in one-to-one correspondence with the sequence of distributions of trees $T_{[n]}$, $n \geq 1$, associated with some homogeneous fragmentation of $\mathbb{N}$.

Before we tackle this (in Proposition 3), we need some more notation. Let $\mathbf{s} = (s_j, j \in \mathbb{N}) \in \mathbb{R}_+^\mathbb{N}$ have total sum $\sum_{j \in \mathbb{N}} s_j \leq 1$. By setting $s_0 = 1 - \sum_{j \in \mathbb{N}} s_j$, we define a probability mass function $(s_j)_{j \geq 0}$ on $\mathbb{N} \cup \{0\}$. Independent random variables $(I_i, i \geq 1)$ with probability mass function $(s_j)_{j \geq 0}$ can be interpreted as an urn scheme, with urns labeled by $\mathbb{N}$ and a "dustbin" with label 0.

As shown by Bertoin, (the laws of) standard homogeneous fragmentations of $\mathbb{N}$ are in one-to-one correspondence with $\sigma$-finite measures $\kappa$ on $\mathcal{P}$ that satisfy

$$\kappa(\{\pi \in \mathcal{P} : \pi|_n \neq \mathbf{1}_{[n]}\}) < \infty \qquad \text{for all } n \geq 1$$

and which informally correspond to the jump measures of the fragmentation processes. We call such measures *dislocation measures* on $\mathcal{P}$. As shown in [8], such measures admit the following, simple, representation. For $\mathbf{s} \in \mathcal{S}^{\downarrow}$, we let $\kappa_{\mathbf{s}}$ be the distribution on $\mathcal{P}$ of the random variable $\Pi$ obtained by Kingman's paintbox construction: let $I_1, I_2 \ldots$ be i.i.d. with law $(s_j)_{j \geq 0}$ and let $i, j$ be in the same block of $\Pi$ if and only if $i = j$ or $I_i = I_j > 0$. Then for every dislocation measure $\kappa$ on $\mathcal{P}$, there exists $c \geq 0$ and a measure $\nu$ on $\mathcal{S}^{\downarrow}$ satisfying (1) such that

$$(8) \qquad \kappa(d\pi) = \int_{\mathcal{S}^{\downarrow}} \kappa_{\mathbf{s}}(d\pi)\nu(d\mathbf{s}) + c\sum_{i=1}^{\infty} \delta_{\epsilon_i}(d\pi),$$

where $\epsilon_i$ is the partition of $\mathbb{N}$ into two blocks $\{i\}$ and $\mathbb{N} \setminus \{i\}$.

2.4. *Characterization of sampling consistent Markov branching models.* We are now almost ready to give the proof of Theorem 1. For any distribution $q_n$ on partitions of the integer $n$ (splitting rule), and for $B$ with $n$ elements, we define the associated *exchangeable splitting rule* on $\mathcal{P}_B \setminus \{\mathbf{1}_B\}$, which is the probability distribution on $\mathcal{P}_B$ defined by

$$(9) \qquad \bar{q}_B(\pi) = \binom{n}{k_1, \ldots, k_r}^{-1} \left(\prod_{i=1}^{n} m_i!\right) q_n(k_1, \ldots, k_r)$$

whenever $\pi$ is a partition of $B$ with $r$ nonempty blocks of sizes $k_1 \geq \cdots \geq k_r$, block size $i$ appearing with frequency $m_i$, $1 \leq i \leq n$. Informally, this is what we obtain when uniformly choosing a partition of $\mathcal{P}_B$ that is compatible with a partition of $n$ that has been sampled according to $q_n$. It is elementary that a random partition with law $\bar{q}_B$ is exchangeable.

Also, it is plain that the laws $(P_n, n \geq 1)$ on labeled trees associated with (not necessarily sampling consistent) splitting rules $(q_n, n \geq 2)$ can also be described as follows. Define $P_{\{i\}}$ to be the Dirac mass on the only element of $\mathbb{T}_{\{i\}}$, in the notation of Section 2.2. Then, recursively, define $P_B$ as the law of $\langle \mathbf{t}_1, \ldots, \mathbf{t}_r \rangle$, where $\pi$ is taken at random according to $\bar{q}_B$ and, given $\pi = (\pi_1, \ldots, \pi_r, \varnothing, \ldots)$ with $\pi_r \neq \varnothing$, $\mathbf{t}_1, \ldots, \mathbf{t}_r$ are picked independently in $\mathbb{T}_{\pi_1}, \ldots, \mathbb{T}_{\pi_r}$ with respective laws $P_{\pi_1}, \ldots, P_{\pi_r}$. Then $P_n = P_{[n]}$. Moreover, $(q_n, n \geq 2)$ is sampling consistent if and only if the image distribution of $P_{[n+1]}$ by the operation that removes the leaf with label $n+1$ is $P_{[n]}$.

PROPOSITION 3. *Sampling consistent splitting rules $(q_n, n \geq 2)$ are in one-to-one correspondence with dislocation measures $\kappa$ on $\mathcal{P}_\mathbb{N}$ of homogeneous fragmentation processes (modulo constant multiples).*

*More precisely, for any $(q_n)_{n \geq 2}$, the formulas $\lambda_2 = 1$,*

$$\lambda_{n+1} = \frac{\lambda_n}{1 - \bar{q}_{[n+1]}(\{1,\ldots,n\},\{n+1\})} \tag{10}$$

*and*

$$\kappa(\{\Gamma \in \mathcal{P} : \Gamma|_n = \pi\}) = \lambda_n \bar{q}_{[n]}(\pi), \tag{11}$$

*for all $\pi \in \mathcal{P}_{[n]} \setminus \{\mathbf{1}_{[n]}\}$, define a dislocation measure $\kappa$ on $\mathcal{P}$.*

*Conversely, for any dislocation measure $\kappa$,*

$$\bar{q}_B(\pi) = \frac{\kappa(\{\Gamma \in \mathcal{P} : \Gamma|_B = \pi\})}{\kappa(\{\Gamma \in \mathcal{P} : \Gamma|_B \neq \mathbf{1}_B\})}, \qquad \pi \in \mathcal{P}_B \setminus \{\mathbf{1}_B\} \text{ and } B \subset_f \mathbb{N}, \tag{12}$$

*defines an exchangeable sampling consistent splitting rule.*

*Moreover, if $\Pi$ is a homogeneous fragmentation process with dislocation measure $\kappa$, then the sequence of distributions of the discrete fragmentation trees $T_{[n]}$, $n \geq 1$, is exactly $(P_n, n \geq 1)$, as associated with $(q_n, n \geq 2)$.*

PROOF. Let $(\Pi(t))_{t \geq 0}$ be a homogeneous $\mathcal{P}$-valued fragmentation process with dislocation measure $\kappa$. For $B \subset_f \mathbb{N}$, let $\bar{q}_B$ be the distribution of $\pi = \Pi|_B(D_B)$, where $D_B = \inf\{t \geq 0 : \Pi|_B(t) \neq \mathbf{1}_B\}$. It is plain that the $\bar{q}_B$ are exchangeable. Thus, they specify partition-valued splitting rules. We denote the associated "unlabeled" splitting rules (i.e., on partitions of $n$) by $q_n$, $n \geq 2$.

By the strong Markov property ([8]) applied at time $D_B$, given $\pi = (\pi_1, \ldots, \pi_r, \varnothing, \ldots)$ with $\pi_r \neq \varnothing$, the processes $(\Pi|_{\pi_i}(D_B + t), t \geq 0)$ for $1 \leq i \leq r$ are independent and, respectively, have the same law as $\Pi|_{\pi_i}, 1 \leq i \leq r$. From this, we see that the discrete fragmentation tree $T_B = \mathbf{t}_{\Pi|_B}$ has distribution $P_B$ associated with the splitting rules $\bar{q}_B$. Sampling consistency for the splitting rules $q_n, n \geq 2$, follows immediately from the property that $T_{[n]}$ is obtained from $T_{[n+1]}$ by deletion of the leaf with label $n+1$. It is argued in Bertoin [8] that $\bar{q}_B$ is indeed given by the formula (12) in the case $B = [n]$ and the general case follows by exchangeability.

Conversely, a sampling consistent system of Markov branching models allows us to consider a strongly sampling consistent system of trees $T_n \sim P_n$, $n \geq 1$. Note that $T_n$ and $T_{n+1}$ are related in one of two possible ways: with probability $p_{n+1} := \bar{q}_{[n+1]}(\{1,\ldots,n\},\{n+1\})$, the branch point adjacent to the root in $T_{n+1}$ splits into $\{1,\ldots,n\}$ and $\{n+1\}$ and has $T_n$ as a subtree; with probability $1 - p_{n+1}$, the branch point closest to the root in $T_{n+1}$ can be identified with the branch point closest to the root of $T_n$. Necessarily, if $P_n, n \geq 0$, can be obtained from some homogeneous fragmentation process $\Pi$, then the holding rates $\lambda_n = \mathbb{E}[D_{[n]}]^{-1}, n \geq 1$, for the state $\mathbf{1}_{[n]}$ of the Markov process $\Pi|_n$ should thus satisfy

$$\mathbb{P}(D_{[n]} = D_{[n+1]}) = 1 - p_{n+1} \quad \text{and, on } \{D_{[n]} \neq D_{[n+1]}\},$$

$$D_{[n]} \sim D_{[n+1]} + \tilde{D}_n,$$

where $\tilde{D}_n$ is independent of $D_{[n+1]}$ and exponential with rate $\lambda_n$. Taking expectations, we get

$$\lambda_n^{-1} = (1 - p_{n+1})\lambda_{n+1}^{-1} + p_{n+1}(\lambda_{n+1}^{-1} + \lambda_n^{-1}) \quad \Rightarrow \quad \lambda_n = (1 - p_{n+1})\lambda_{n+1}.$$

If we arbitrarily put $\lambda_2 = 1$, this determines $(\lambda_n)_{n \geq 3}$ from $(q_n)_{n \geq 2}$. Furthermore, by the same reasoning, we get, for all $\pi \in \mathcal{P}_{[n]}$,

$$\bar{q}_{[n]}(\pi) = \bar{q}_{[n+1]}(\{\Gamma \in \mathcal{P}_{[n+1]} : \Gamma|_n = \pi\})$$
(13)
$$+ \bar{q}_{[n+1]}(\{1, \ldots, n\}, \{n+1\})\bar{q}_{[n]}(\pi),$$

that is, after rearrangement and multiplication by $\lambda_{n+1}$,

$$\lambda_n \bar{q}_{[n]}(\pi) = \lambda_{n+1} \bar{q}_{[n+1]}(\{\Gamma \in \mathcal{P}_{[n+1]} : \Gamma|_n = \pi\})$$

so that we can define $\kappa$ consistently by (11) as a $\sigma$-finite measure on $\mathcal{P}_\mathbb{N}$.

Finally, for the uniqueness, note that the choice of $\lambda_2$ was arbitrary and any other choice $\lambda_2 \in (0, \infty)$ leads to a constant multiple of $\kappa$, that is, a linear time change of an associated fragmentation process. It is easily checked that if $\kappa$ is defined by (10) and (11) for any $\lambda_2 \in (0, \infty)$, then (12) holds for $B = [n]$ and then for $B \subset [n]$; and if $\bar{q}_{[n]}$ is defined by (12), then (10) and (11) hold with $\lambda_n = \kappa(\{\Gamma \in \mathcal{P}_\mathbb{N} : \Gamma_n \neq \mathbf{1}_{[n]}\})$. $\square$

The consistency equation (13) can be written in terms of $q_n$ as

(14)
$$q_n(k_1, \ldots, k_r) = \sum_{j=1}^r \frac{(k_j + 1)(m_{k_j+1} + 1)}{(n+1)m_{k_j}} q_{n+1}((k_1, \ldots, k_j + 1, \ldots, k_r)^\downarrow)$$
$$+ \frac{m_1 + 1}{n + 1} q_{n+1}(k_1, \ldots, k_r, 1)$$
$$+ \frac{1}{n+1} q_{n+1}(n, 1) q_n(k_1, \ldots, k_n),$$

which is structurally similar to but not the same as, the backward recursions for the rows of the decrement matrix associated with coalescents with simultaneous multiple collisions; see [13]. The binary special case was already obtained by Ford [21], Proposition 41, and can be compared with coalescents with no simultaneous but multiple collisions, as in [14]. See also [24, 25] for similar recursions in the context of regenerative composition and partition structures.

PROOF OF THEOREM 1. The fact that all sampling consistent splitting rules are of the stated form is now a simple exercise using (12), (8) and (9). The theorem will be proven if we show that $c, \nu$ can be recovered from the dislocation measure $\kappa$ associated, up to a constant multiple, with a sampling consistent splitting rule as in Proposition 3. Obviously, $c = \kappa(\{\epsilon_1\})$, so we can assume $c = 0$ in (8). We then use Kingman's paintbox construction to obtain that $\kappa$-almost every $\pi$ has an asymptotic frequency, and that the restriction of $\kappa$ to $A_\varepsilon := \{\pi \in \mathcal{P} : \max_i |\pi_i| < 1 - \varepsilon\}$ is finite with total mass $m_\varepsilon = \nu(\{\mathbf{s} \in \mathcal{S}^\downarrow : s_1 < 1 - \varepsilon\})$. Then the probability measure $\nu(\cdot \cap \{\mathbf{s} : s_1 < 1 - \varepsilon\})/m_\varepsilon$ is just the distribution of $|\pi|^\downarrow$ under $\kappa(\cdot \cap A_\varepsilon)/m_\varepsilon$ so that $\nu$ is recovered from $\kappa$. □

The binary special case is worth discussing separately. It is characterized by those dislocation measures that have the property

$$\kappa(\{\pi = (\pi_1, \pi_2, \ldots) \in \mathcal{P}_\mathbb{N} : \pi_1 \cup \pi_2 \neq \mathbb{N}\}) = 0.$$

Writing $\kappa = \kappa_0 + c \sum_{i \geq 1} \delta_{(\{i\}, \mathbb{N} \setminus \{i\})}$ (for the highest $c$ such that $\kappa_0$ is a nonnegative measure) and using the one-to-one correspondence of dislocation measures on $\mathcal{P}_\mathbb{N}$ and pairs of erosion coefficient $c \geq 0$ and dislocation measure $\nu$ on $\mathcal{S}^\downarrow$, these correspond to $(c, \nu)$ with

(15) $\qquad c \geq 0 \quad \text{and} \quad \nu(\{(s_i)_{i \geq 1} \in \mathcal{S}^\downarrow : s_1 + s_2 < 1\}) = 0.$

The presentation is nicest for a symmetrized setting. We define $\tilde{\nu}(A) = \frac{1}{2}(\nu(s_1 \in A) + \nu(s_2 \in A))$ for Borel sets $A \subseteq [0, 1]$.

COROLLARY 4. *Sampling consistent binary splitting rules $q_n, n \geq 2$, are in one-to-one correspondence (modulo constant multiples) with pairs $(c, \nu)$ satisfying* (15).

*Specifically, for any $(c, \nu)$,*

$$\tilde{q}_n(k) = \frac{1}{Z_n}\left(\binom{n}{k} \int_{(0,1)} x^k (1-x)^{n-k} \tilde{\nu}(dx) + nc\mathbf{1}_{\{k=1\}}\right),$$

$1 \leq k \leq n-1$, *where $Z_n = \int_{(0,1)} (1 - x^n - (1-x)^n) \tilde{\nu}(dx) + nc$ is the normalization constant, induces a sampling consistent splitting rule by $q_n(n-k, k) = \tilde{q}_n(k) + \tilde{q}_n(n-k)$, $1 \leq k < n/2$, $q_n(n/2, n/2) = \tilde{q}_n(n/2)$, $n$ even.*

The symmetric splitting rules $\tilde{q}_n(k)$ [for $c = 0$ and $\tilde{\nu}(dx) = f(x)\,dx$ absolutely continuous] give Aldous's (planar) Markov branching models and Corollary 4 shows that, essentially, Aldous had found all binary exchangeable sampling consistent Markov branching models without erosion and expressed them in terms of (a density of) a binary dislocation measure.

Vice versa, we can calculate $c$ and $\nu$ from $n^{-1}Z_n\tilde{q}_n(1) \to c$ and

$$(16) \quad Z_n \sum_{an \leq k \leq bn} \tilde{q}_n(k) = \sum_{an \leq k \leq bn} \int_0^1 \binom{n}{k} x^k(1-x)^{n-k}\tilde{\nu}(dx) \to \tilde{\nu}([a,b])$$

for any continuity points $0 < a < b < 1$ for $\tilde{\nu}$, provided $Z_n$ (or another normalization sequence $\tilde{Z}_n \sim Z_n$) can be calculated. The proof, which is easily done using de Finetti's representation for exchangeable sequences of 0's and 1's, is left as an exercise to the reader.

**3. Self-similar fragmentations and continuum trees.** In this section, we set the bases needed to prove convergence of discrete fragmentation trees to some continuum random trees that are naturally related to the so-called self-similar fragmentations [9, 10].

3.1. *Self-similar fragmentations.* A nice feature of exchangeable partitions in the case $B = \mathbb{N}$ is that Kingman's theory [31] entails that the blocks of such a partition $\Pi$ admit *asymptotic frequencies* almost surely, namely

$$|\Pi_i| := \lim_{n\to\infty} \frac{\#\Pi_i \cap [n]}{n}.$$

We let $|\Pi| = (|\Pi_i|, i \geq 1)$ and $|\Pi|^{\downarrow}$ be the random element of $\mathcal{S}^{\downarrow}$ obtained from $|\Pi|$ by ranking its terms in decreasing order.

Let $(\Pi(t), t \geq 0)$ be an exchangeable càdlàg (right-continuous with left limits) $\mathcal{P}$-valued stochastic process such that $\Pi(0) = \mathbf{1}_{\mathbb{N}}$, and $|\Pi(t)|$ exists for every $t \geq 0$, a.s. Suppose, also, that the process of sizes of the block containing $i$, $(|\Pi_{(i)}(t)|, t \geq 0)$, is right-continuous for every $i \in \mathbb{N}$ a.s.

DEFINITION 4. The process $(\Pi(t), t \geq 0)$ is a $\mathcal{P}$-valued self-similar fragmentation process with index $a \in \mathbb{R}$ if, given $\Pi(t) = \pi$, the random variable $\Pi(t+s)$ has same law as the random partition whose blocks are those of $\pi_i \cap \Pi^{(i)}(|\pi_i|^a s), i \geq 1$, where $(\Pi^{(i)}, i \geq 1)$ is a sequence of i.i.d. copies of $(\Pi(t), t \geq 0)$.

When $a = 0$, we recover the definition of standard homogeneous fragmentations in $\mathcal{P}$. To avoid trivialities, we will only work with nonconstant processes. We notice that if $(\Pi(t), t \geq 0)$ is a self-similar $\mathcal{P}$-valued fragmentation, then $(|\Pi(t)|^{\downarrow}, t \geq 0)$ is a self-similar fragmentation with values in $\mathcal{S}^{\downarrow}$, as defined in the Introduction (and any $\mathcal{S}^{\downarrow}$-valued fragmentation can be represented in this form; see [7, 9]). Bertoin has shown in [9] that $\mathcal{P}$-valued self-similar fragmentations are characterized by a triple $(a, c, \nu)$, where $c \geq 0$ and $\nu$ is a dislocation measure (1) on $\mathcal{S}^{\downarrow}$. Hereafter, we will only be interested in the cases where $c = 0$ and $\nu$ is conservative (4) (no sudden loss of mass and no erosion). We call $(a, \nu)$ the *characteristic pair* of such self-similar fragmentations.

There is a useful way to relate self-similar fragmentations to homogeneous fragmentations, which is as follows.

LEMMA 5 ([9]). *Let $(\Pi^0(t), t \geq 0)$ be a standard homogeneous fragmentation with dislocation measure $\nu$ and let $a \in \mathbb{R}$. We then define a sequence of time changes,*

$$(17) \qquad \eta_{(i)}(t) = \inf\left\{u \geq 0 : \int_0^u |\Pi^0_{(i)}(w)|^{-a}\, dw > t\right\}, \qquad t \geq 0,\ i \geq 1.$$

*Let $\Pi^a(t)$ be the element of $\mathcal{P}$ whose blocks are those of the partitions $\Pi^0_{(i)}(\eta_{(i)}(t))$, $i \geq 1$. Then:*

   (i) *the process $(\Pi^a(t), t \geq 0)$ is a self-similar fragmentation with characteristic pair $(a, \nu)$;*

   (ii) *for the size $|\Pi^0_{(i)}(t)|$ of the block containing $i$, the process $\xi_{(i)}(t) = -\log|\Pi^0_{(i)}(t)|, t \geq 0$ is a pure-jump subordinator with Lévy measure*

$$(18) \qquad \Lambda(dx) = e^{-x} \sum_{i \geq 1} \nu((s_j)_{j \geq 1} \in \mathcal{S}^{\downarrow} : -\log s_i \in dx).$$

*Thus, $|\Pi^a_{(i)}(t)| = \exp(-\xi_{(i)}(\eta_{(i)}(t)))$, where*

$$(19) \qquad \eta_{(i)}(t) = \inf\left\{u \geq 0 : \int_0^u e^{a\xi_{(i)}(w)}\, dw > t\right\}, \qquad t \geq 0.$$

We refer to [9] for the proof of this result. Note that because the partitions $\Pi^0(t)$ are refining as $t$ increases, if two of the blocks of the partitions $\Pi^0_{(i)}(\eta_{(i)}(t))$, $i \geq 1$, have a common element, then they are equal and the definition of $\Pi^a(t)$ makes sense.

3.2. *Trees with edge lengths.* Let $(\Pi(t), t \geq 0)$ be a self-similar fragmentation process with index $a$. We may then construct a family of random trees $T_B$ indexed by $B \subset_f \mathbb{N}$, defined by $T_B = \mathbf{t}_{\Pi|_B}$, the fragmentation tree associated with the restrictions of $(\Pi(t), t \geq 0)$ to $B$ (see Definition 1). The time-change construction of Lemma 5 provides a coupling of all self-similar fragmentations with same dislocation measure and different indices $a \in \mathbb{R}$, all with the same $T_B$. The only difference lies in the times at which splits occur, which do not appear in $T_B$. These times provide extra information on the tree associated with a fragmentation process, which we can interpret as *edge lengths* associated with the fragmentation tree $T_B$ for a particular index $a \in \mathbb{R}$.

A (rooted labeled) *tree with edge lengths* is a pair $\vartheta = (\mathbf{t}, (\eta_e, e \in E(\mathbf{t})))$, where $\mathbf{t} \in \mathbb{T}_B$ for some $B \subset_f \mathbb{N}$, $E(\mathbf{t})$ is the set of edges of $\mathbf{t}$ and $(\eta_e, e \in E(\mathbf{t})) \in (0, \infty)^{E(\mathbf{t})}$ are positive marks, interpreted as the lengths of the associated edges. The tree $\mathbf{t}$ is called the *shape* and we let $\Theta_B$ be the set of trees with edge lengths whose shape is in $\mathbb{T}_B$.

Let $(\pi(t), t \geq 0)$ be a refining process with values in $\mathcal{P}$, starting at $\mathbf{1}_\mathbb{N}$. Assume, further, that

$$D_i := D_{\{i\}} = \inf\{s \geq 0 : \{i\} \in \pi(s)\} < \infty \qquad \text{for all } i \geq 1,$$

so that in particular, $\pi|_B(t) = \mathbf{0}_B$ for some finite $t$. Recall that a vertex $v$ of any $\mathbf{t} \in \mathbb{T}_B$ is naturally identified with the set $B_v$ of labels of the leaves that descend from $v$. We are going to make this identification in the sequel.

For $B \subset_f \mathbb{N}$, we let $\theta_{\pi|_B} \in \Theta_B$ be the tree with edge lengths whose shape is $\mathbf{t}_{\pi|_B}$ and whose edge lengths are $\eta_e = D_v - D_{\neg v}$ whenever $e \in E(\mathbf{t}_{\pi|_B})$ is the edge linking a nonroot vertex $v$ and its parent $\neg v$. Notice that whereas $D_v = \inf\{t \geq 0 : \Pi|_{B_v}(t) \neq \mathbf{1}_{B_v}\}$ for a nonleaf vertex $v$ only depends on $\Pi|_{B_v}$, $D_i$ is defined differently and depends on the whole process $(\pi(t), t \geq 0)$ rather than its restriction to $\{i\}$ or $B$.

Now, suppose that $(\Pi(t), t \geq 0)$ is a self-similar fragmentation with dislocation measure $\nu$ and index $a < 0$. By [9], it holds that $0 < D_i < \infty$ a.s. for every $i$ and, in fact, $\sup_{i \geq 1} D_i < \infty$ in that case. Therefore, $\mathcal{R}_B = \theta_{\Pi|_B}$ is well defined. This tree was called $\mathcal{R}(B)$ in [29], Section 2.3, where it was constructed slightly differently. We conclude this section by establishing the link between the two presentations.

If $\vartheta \in \Theta_B$ has a root-edge $e$ with length $\eta_e$ and if $x < \eta_e$, then we let $\vartheta - x$ be the element of $\Theta_B$ with the same shape and edge lengths, except for the root-edge, which is assigned length $\eta_e - x$. If $\vartheta_1, \ldots, \vartheta_r$ are elements of $\Theta_{B_1}, \ldots, \Theta_{B_r}$ with shapes $\mathbf{t}_1, \ldots, \mathbf{t}_r$ for pairwise disjoint nonempty $B_i \subset_f \mathbb{N}$, and if $x > 0$, we let

$$\langle \vartheta_1, \ldots, \vartheta_r \rangle_x$$

be the element of $\Theta_{\cup_i B_i}$ whose shape is $\langle \mathbf{t}_1, \ldots, \mathbf{t}_r \rangle$, whose root-edge length is $x$ and whose other edge lengths are inherited from those of $\vartheta_1, \ldots, \vartheta_r$ in the natural way.

The trees $\mathcal{R}_B$ can be recursively described in the following way. Let $\mathcal{R}_{\{i\}}$ have as shape the only element of $\mathbb{T}_{\{i\}}$ and (single) edge length equal to $D_{\{i\}}$. Then let

$$\mathcal{R}_B = \langle \mathcal{R}_{B_1} - D_B, \ldots, \mathcal{R}_{B_r} - D_B \rangle_{D_B},$$

where $B_1, \ldots, B_r$ are the nonempty blocks of the partition of $B$ induced by $\Pi(D_B)$. This is the definition of [29], Section 2.3.

3.3. *Continuum trees and fragmentation processes.*

3.3.1. $\mathbb{R}$-*trees.* We now introduce the continuous version of trees that is needed to deal with continuum random trees, following [4, 19].

An $\mathbb{R}$-tree $(\tau, d)$ is a complete separable metric space such that for every $x, y \in \tau$:

1. there is an isometry $\varphi_{x,y} \colon [0, d(x,y)] \to \tau$ such that $\varphi_{x,y}(0) = x$ and $\varphi_{x,y}(d(x,y)) = y$;

2. for every injective path $c \colon [0,1] \to \tau$ with $c(0) = x$, $c(1) = y$, one has $c([0,1]) = \varphi_{x,y}([0, d(x,y)])$.

In other words, there exists a geodesic in $\tau$ linking any two points and this geodesic is the only simple path linking these points (up to reparameterization). We usually denote by $[[x, y]]$ the range of $\varphi_{x,y}$. This is indeed a continuous analog of the graph-theoretic definition of a tree as a connected graph with no cycle. The $\mathbb{R}$-trees we will be considering are also *rooted*, that is, they have a distinguished element which we denote by $\rho$.

We say that two rooted $\mathbb{R}$-trees $(\tau, \rho, d), (\tau', \rho', d')$ are *equivalent* if there exists an isometry from $\tau$ onto $\tau'$ that sends the $\rho$ to $\rho'$. We denote by $\Theta$ the set of equivalence classes of compact rooted $\mathbb{R}$-trees. It has been shown in [19] that $\Theta$ is a Polish space when endowed with the so-called *pointed Gromov–Hausdorff distance* $d_{\mathrm{GH}}$, where, by definition, the distance $d_{\mathrm{GH}}((\tau, \rho), (\tau', \rho'))$ is equal to the infimum of the quantities

$$\delta(r, r') \vee \delta_{\mathcal{H}}(T, T'),$$

where $(T, r), (T', r')$ are isometric embeddings of $(\tau, \rho), (\tau', \rho')$ into a common metric space $(M, \delta)$ and $\delta_{\mathcal{H}}$ is the Hausdorff distance between compact subsets of $(M, \delta)$. It is elementary that this distance does not depend on particular choices in the equivalence classes of $(\tau, \rho)$ and $(\tau', \rho')$. We endow $\Theta$ with the associated Borel $\sigma$-algebra. In the sequel, by a slight abuse of notation, we will still call elements of $\Theta$ rooted $\mathbb{R}$-trees, and we will denote them by $\tau$, omitting mention of the root and the distance $d$. Also, by a probability measure on an element $\tau \in \Theta$, we will mean an equivalence class of a 4-tuple $(\tau, \rho, d, \mu)$, where we call $(\tau, \rho, d, \mu)$ and $(\tau', \rho', d', \mu')$ equivalent if there exists an isometry from $(\tau, \rho, d)$ to $(\tau', \rho', d')$ such that $\mu'$ is the push-forward of $\mu$.

If $\tau \in \Theta$, then and for $x \in \tau$, we call the quantity $d(\rho, x)$ the *height* of $x$. If $x, y \in \tau$, we say that $x$ is an ancestor of $y$ whenever $x \in [[\rho, y]]$. We let $x \wedge y \in \tau$ be the unique element of $\tau$ such that $[[\rho, x]] \cap [[\rho, y]] = [[\rho, x \wedge y]]$ and call it the highest common ancestor of $x$ and $y$ in $\tau$. For $x \in \tau$, we denote by $\tau_x$ the set of $y \in \tau$ such that $x$ is an ancestor of $y$. The set $\tau_x$, endowed with the restriction of the distance $d$ and rooted at $x$, is in turn a rooted $\mathbb{R}$-tree called the *fringe subtree* of $\tau$ rooted at $x$.

We say that $x \in \tau$, $x \neq \rho$, in a rooted $\mathbb{R}$-tree is a *leaf* if its removal does not disconnect $\tau$ and we let $\mathcal{L}(\tau)$ denote the set of leaves of $\tau$. A *branch point* is an element of $\tau$ of the form $x \wedge y$, where $x$ is not an ancestor of $y$, nor vice versa. It is also characterized by the fact that the removal of a branch point disconnects the $\mathbb{R}$-tree into three or more components (two or more for the root). We let $\mathcal{B}(\tau)$ denote the set of branch points of $\tau$.

3.3.2. *Relation with trees with edge lengths, reduced trees.* There is a natural connection between the trees with edge lengths with shape in $\mathbb{T}_n^\circ$ (resp. $\mathbb{T}_B$) of the previous sections and rooted $\mathbb{R}$-trees with $n$ leaves (resp. $\#B$ leaves labeled by $B$) and where the root is not a branch point. If $\tau$ is a rooted $\mathbb{R}$-tree with $\rho \notin \mathcal{B}(\tau)$ and exactly $n$ leaves labeled $L_1, \ldots, L_n$, then we consider the graph whose vertices are the set $V = \{\rho\} \cup \mathcal{L}(\tau) \cup \mathcal{B}(\tau)$ and such that two vertices $x, y$ are connected by an edge if and only if $[[x, y]] \cap V = \{x, y\}$. The resulting graph is a tree which is naturally rooted at $\rho$ and the edge connecting $x$ and $y$ naturally inherits the length $d(x, y) = |d(\rho, x) - d(\rho, y)|$. This construction can be reversed, associating an $\mathbb{R}$-tree with a tree with edge lengths, for example, by means of Aldous's *sequential construction* [4], page 252.

Also, if **t** is an element of $\mathbb{T}_n^\circ$ or $\mathbb{T}_B$, one naturally puts edge lengths equal to 1 on each edge and considers **t** as an $\mathbb{R}$-tree as well, the restriction of the distance of this $\mathbb{R}$-tree to the set of branch points, leaves and the root being the usual combinatorial distance on the vertices of **t**.

For $\tau$ a rooted $\mathbb{R}$-tree and $x_1, x_2, \ldots, x_n \in \tau$, we let

$$R(\tau, x_1, \ldots, x_n) = \bigcup_{i=1}^{n} [[\rho, x_i]]$$

be the *reduced* subtree associated with $\tau, x_1, \ldots, x_n$. It is elementary that $R(\tau, x_1, \ldots, x_n)$ is, in turn, an $\mathbb{R}$-tree, which is naturally rooted at $\rho$ and whose leaves are included in $\{x_1, \ldots, x_n\}$ (it might be that $x_i$ is not a leaf of the reduced tree whenever $x_i$ is an ancestor of $x_j$ for some $j \neq i$, but note that this never happens if $x_1, \ldots, x_n$ are distinct leaves of $\tau$). By the discussion of the previous paragraph, if $\tau$ is such that $\rho \notin \mathcal{B}(\tau)$ and if $x_1, \ldots, x_n$ are leaves of $\tau$, then $R(\tau, x_1, \ldots, x_n)$ can also be considered as a tree with edge lengths, whose shape is in $\mathbb{T}_n$, since the leaves inherit a natural labelling from that of $x_1, \ldots, x_n$.

3.3.3. *Continuum random trees and fragmentation trees.* The fragmentation trees introduced in [29] are yet another way to consider self-similar fragmentation processes whose characteristic pair $(-\gamma, \nu)$ satisfies $\gamma > 0$. In order to introduce them, we first need to give some definitions and results on continuum trees, following [4].

We say that a pair $(\tau, \mu)$ is a *continuum tree* if $\tau \in \Theta$ and $\mu$ is a probability measure on $\tau$ such that:

1. $\mu$ is supported by the set $\mathcal{L}(\tau)$;
2. $\mu$ has no atom;
3. for every $x \in \tau \setminus \mathcal{L}(\tau)$, $\mu(\tau_x) > 0$.

Notice that continuum trees automatically satisfy a number of properties. For example, the set $\mathcal{L}(\tau)$ must be uncountable (by 1 and 2) and cannot have isolated points (by 2 and 3).

A *continuum random tree* (CRT) is a "random variable" whose values are continuum trees, defined on some probability space $(\Omega, \mathcal{A}, \mathbb{P})$. To formalize this, we should endow the set of continuum trees with a $\sigma$-algebra. A natural possibility would be to use Evans' and Winter's separable and complete metric structure [20] on the space of "weighted $\mathbb{R}$-trees," although we would have to incorporate the fact that our trees are rooted. Another, probably even more natural, approach would be to use the Gromov-weak topology on the set of metric measure spaces introduced in the recent work of Greven, Pfaffelhuber and Winter [27].

However, for technical simplicity, in this paper, we prefer to follow Aldous [4] and use the space $l_1 = l_1(\mathbb{N})$ as a base space for defining our CRTs. Namely, we endow the set of compact subsets of $l_1$ with the Hausdorff metric, and the set of probability measures on $l_1$ with any metric inducing the topology of weak convergence, so that the set of pairs $(T, \mu)$, where $T$ is a rooted $\mathbb{R}$-tree embedded as a subset of $l_1$ and $\mu$ is a measure on $T$, is endowed with the product $\sigma$-algebra. Convergence for the Hausdorff metric for subsets of $l_1$ is stronger than convergence of the associated equivalence classes for the Gromov–Hausdorff topology. In the sequel, we always keep in mind that the usual probability "operations" such as conditioning, for example, with respect to $\mu$, and then sampling i.i.d. random variables with law $\mu$, are done by using this $l_1$-embedded measurable representative before taking isometric equivalence classes. In this sense, given $(\mathcal{T}, \mu)$, we will call an i.i.d. sequence $L_1, L_2, \ldots$ with common law $\mu$ an *exchangeable sequence with directing law $\mu$*.

For $a > 0$ and $(\tau, \rho(\tau), d) \in \Theta$, we denote by $a\tau$ the element $(\tau, \rho(\tau), ad)$ obtained by scaling distances by a factor $a$.

For $(\tau, \mu)$ a continuum tree, we let $C_t^1, C_t^2, \ldots$ be the connected components of the open set $\{x \in \tau : d(x, \rho(\tau)) > t\}$, ranked so that $\mu(C_t^1) \geq \mu(C_t^2) \geq \cdots$. We then let $\sigma_t^i$ be the element of $\tau$ at height $t$ such that $C_t^i \subset \tau_{\sigma_t^i}$. Then $\tau_t^i = C_t^i \cup \{\sigma_t^i\}$ is a compact $\mathbb{R}$-tree which we root at $\sigma_t^i$. Notice that $\tau_t^i$ is equal to $\tau_{\sigma_t^i}$ unless $\sigma_t^i \in \mathcal{B}(\tau)$.

DEFINITION 5. *A self-similar tree with index $-\gamma < 0$ is a continuum random tree $(\mathcal{T}, \mu)$ such that for every $t \geq 0$, given $(\mu(\mathcal{T}_t^i), i \geq 1)$, the continuum random trees*

$$\left(\mu(\mathcal{T}_t^1)^{-\gamma}\mathcal{T}_t^1, \frac{\mu(\cdot \cap \mathcal{T}_t^1)}{\mu(\mathcal{T}_t^1)}\right), \left(\mu(\mathcal{T}_t^2)^{-\gamma}\mathcal{T}_t^2, \frac{\mu(\cdot \cap \mathcal{T}_t^2)}{\mu(\mathcal{T}_t^2)}\right), \ldots$$

are i.i.d. copies of $(\mathcal{T}, \mu)$.

Again, the last sentence means that there exist i.i.d. copies of a representative of $(\mathcal{T}, \mu)$ embedded in $l_1$, whose isometry classes are those of $(\mu(\mathcal{T}_t^1)^{-\gamma}\mathcal{T}_t^1, \mu(\cdot \cap \mathcal{T}_t^1)/\mu(\mathcal{T}_t^1)), \ldots$.

As was shown in [29], Theorem 1, Proposition 1, (laws of) self-similar continuum random trees with index $-\gamma < 0$ are in one-to-one correspondence with (laws of) self-similar fragmentation processes with index $-\gamma$, no erosion and no sudden loss of mass. We briefly describe how the two objects are related.

PROPOSITION 6. *Let $(-\gamma, \nu)$ be a characteristic pair with $\gamma > 0$. There then exists a (unique) self-similar CRT $(\mathcal{T}, \mu)$ with index $-\gamma$ such that the following holds. Given $(\mathcal{T}, \mu)$, let $L_1, L_2, \ldots$ be an exchangeable sequence with directing law $\mu$. For every $t \geq 0$, we let $\Pi(t)$ be the random element of $\mathcal{P}$ such that $i$ and $j$ are in the same block of $\Pi(t)$ if and only if $d(\rho(\mathcal{T}), L_i \wedge L_j) > t$, that is, if and only if $L_i$ and $L_j$ belong to the same element of $\{\mathcal{T}_t^1, \mathcal{T}_t^2, \ldots\}$. Then:*

(i) *the process $(\Pi(t), t \geq 0)$ is a $\mathcal{P}$-valued self-similar fragmentation with characteristic pair $(-\gamma, \nu)$ and the process $((\mu(\mathcal{T}_t^i), i \geq 1), t \geq 0)$ is equal to the process $(|\Pi(t)|^{\downarrow}, t \geq 0)$, that is, it is an $\mathcal{S}^{\downarrow}$-valued fragmentation process with characteristic pair $(-\gamma, \nu)$;*

(ii) *if $\mathcal{T}_t^{(i)}$ denotes the element of $\{\mathcal{T}_t^1, \mathcal{T}_t^2, \ldots\}$ that contains $L_i$, then the process $(\mu(\mathcal{T}_t^{(i)}), t \geq 0)$ is equal to $(|\Pi_{(i)}(t)|, t \geq 0)$;*

(iii) *the reduced tree $R(\mathcal{T}, L_1, \ldots, L_k)$ is equal to the tree with edge lengths $\mathcal{R}_{[k]}$, as defined in Section 3.2.*

PROOF. (i) It was shown in [29] that there is a unique tree $(\mathcal{T}, \mu)$ such that $((\mu(\mathcal{T}_t^i), i \geq 1), t \geq 0)$ is the $\mathcal{S}^{\downarrow}$-valued fragmentation process with characteristic pair $(-\gamma, \nu)$. The fact that $(\mu(\mathcal{T}_t^i), i \geq 1) = |\Pi(t)|^{\downarrow}$ for every $t$ comes from the fact that $\mu$ is a.s. the limit of the empirical measure on $L_1, L_2, \ldots$. It is easy to show that this process is right-continuous and that $(\Pi(t), t \geq 0)$ is a càdlàg $\mathcal{P}$-valued process, and it follows that $(\Pi(t), t \geq 0)$ is the $\mathcal{P}$-valued fragmentation process with characteristic pair $(-\gamma, \nu)$. (ii) is immediate from the fact that $i$ and $j$ are in the same block of $\Pi(t)$ if and only if $L_j$ is in $\mathcal{T}_t^{(i)}$ and the fact that $\mu$ is a.s. the limit of the empirical measure on the leaves $L_1, \ldots, L_n$ as $n \to \infty$. Finally, (iii) is [29], Lemma 4. □

**4. Asymptotics of discrete fragmentation trees.** We now embark on the proof of Theorem 2. We can obtain a weaker statement of convergence in distribution by using Aldous' theorems [4], Theorem 18, Corollary 19 and Remark 4. With a little more effort, we establish the stronger statement of Theorem 2 that, in fact, there exists a fragmentation tree defined on the given probability space, to which the discrete fragmentation trees converge in probability.

It will often be convenient to *initially* assume the following.

HYPOTHESIS (H). Assume that we are given a probability space supporting $(\mathcal{T}, \mu)$, a fragmentation tree associated with a self-similar fragmentation with characteristic pair $(-\gamma_\nu, \nu)$, where $\nu$ satisfies the assumptions (6), (7). For simplicity, we let $\gamma = \gamma_\nu$. We assume that our probability space also supports $L_1, L_2, \ldots$, an exchangeable sample of leaves with directing measure $\mu$. We let $\mathcal{R}_k = R(\mathcal{T}, L_1, \ldots, L_k)$ and define a self-similar $\mathcal{P}$-valued fragmentation process $(\Pi(t), t \geq 0)$ with index $-\gamma$ by the device explained in Proposition 6. Also, we let $(\Pi^0(t), t \geq 0)$ be the homogeneous fragmentation process obtained from $(\Pi(t), t \geq 0)$ by the time-change transformation of Lemma 5. We denote by $\xi(t) = -\log |\Pi^0_{(1)}(t)|, t \geq 0$ the pure-jump subordinator with Lévy measure (18) associated by Lemma 5.

We let $T_n$ be the discrete fragmentation tree with $n$ leaves associated with $(\Pi(t), t \geq 0)$ [or $(\Pi^0(t), t \geq 0)$], as in Section 2. The tree $T_n$ is then considered as an $\mathbb{R}$-tree by assuming that its edges are segments with length 1, in accordance with the discussion of Section 3.3.2.

To see that Theorem 2 remains true without Hypothesis (H), simply note that a strongly consistent sequence $(T_n^\circ)$ has the same distribution (as a sequence of random variables) as if it were constructed under Hypothesis (H). Since convergence in probability for random variables with values in a complete space can be metrized by a complete distance (see [15], Theorem 9.2.3), we deduce that $(T_n^\circ)$ is a Cauchy sequence for this distance, and thus converges in probability, because the set of compact real trees endowed with the Gromov–Hausdorff distance is complete. The distribution of the limit is identified as that of the fragmentation tree $\mathcal{T}$.

We recall that all trees involved in Aldous' study [4] are subspaces of $l_1$, endowed with the $l_1$-distance, and that convergence of (compact) trees holds with respect to the Hausdorff distance. Using $l_1$-representatives of the trees $T_n$ and $\mathcal{T}$ (which is always possible; see [29]) and applying Aldous' asymptotic results will then lead us to the claimed convergence in the Gromov–Hausdorff sense. More precisely, using Theorem 18, Corollary 19 and Remark 4 in [4], we see that the proof of the convergence in distribution for the Gromov–Hausdorff topology,

$$\frac{T_n}{n^\gamma \ell(n)} \xrightarrow{(d)} \Gamma(1-\gamma)\mathcal{T},$$

amounts to the following:

(i) the leaf-tightness of $(\mathcal{R}_k, k \geq 1)$, that is, $\min_{2 \leq j \leq k} d(L_1, L_j) \xrightarrow{(p)} 0$ as $k \to \infty$;

(ii) the (a.s.) compactness of $\mathcal{T}$;

(iii) the convergence of "finite-dimensional marginals";

(iv) a tightness criterion, which is stated precisely in Proposition 9 below.

We will obtain the stronger convergence in probability under Hypothesis (H) by using the particular coupling of the discrete and continuum fragmentation trees, and establishing an almost sure convergence result in (iii). The tightness estimate of Aldous then provides the uniform bound that is needed to extend convergence of finite-dimensional marginals to Gromov–Hausdorff convergence, at the price of turning the a.s. convergence into convergence in probability.

The two first points are proved in [29], Lemmas 3 and 5. The aim of this section is therefore to prove the latter two: Section 4.1 is devoted to the convergence of finite-dimensional marginals and Section 4.2 to the tightness estimate. Section 4.2 also contains the proof of Theorem 2. Finally, Section 4.3 provides an analog of Theorem 2 for convergence of leaf-height functions.

4.1. *Convergence of finite-dimensional marginals.* The first step is given by the following proposition, which contains the convergence of "finite dimensional marginals" for Theorem 2, but note that we do not need the integrability condition (7).

PROPOSITION 7. *Let $\nu$ be a conservative dislocation measure satisfying the regular variation condition (6), $(T_n)_{n \geq 1}$ an associated strongly sampling consistent family of discrete fragmentation trees defined on any probability space. Then the same probability space also supports $\mathcal{R}_k$ so that*

$$n^{-\gamma} \ell(n)^{-1} R(T_n, [k]) \xrightarrow[n \to \infty]{\text{a.s.}} \Gamma(1 - \gamma) \mathcal{R}_k,$$

*in the Gromov–Hausdorff sense, for all $k \geq 1$.*

We observe that the convergence is in the sense of the Gromov–Hausdorff metric, but in the context of trees with edge lengths, there is a simple sufficient condition: finite trees with edge lengths $\vartheta_n$ converge to another finite tree with edge lengths $\vartheta$ if the shape of $\vartheta_n$ is eventually that of $\vartheta$ and the edge lengths converge pointwise. This condition is almost necessary, but there is a complication when some edge lengths converge to zero and shapes oscillate—this will be irrelevant here.

A key ingredient is provided by the following lemma.

LEMMA 8 ([26]). *Let $\xi = (\xi_t, t \geq 0)$ be a pure-jump subordinator with Lévy measure $\Lambda$ satisfying*

(20) $$\Lambda([x, \infty)) = x^{-\gamma} \ell(1/x), \qquad x \downarrow 0.$$

Let $V_1, V_2, \ldots$ be a sequence of nonnegative random variables which conditionally given $\xi$ are independent and identically distributed with

$$\mathbb{P}(V_i > v) = \exp(-\xi_v), \qquad v \geq 0,$$

and for $s \geq 0$, let

$$K_n(s) := \#\{V_i : 1 \leq i \leq n, V_i \leq s\},$$

which is the number of distinct values among the $V_i : 1 \leq i \leq n$ with $V_i \leq s$. Then

$$(21) \quad \lim_{n \to \infty} \sup_{0 \leq s \leq \infty} \left| \frac{K_n(s)}{n^\gamma \ell(n)} - \Gamma(1-\gamma) \int_0^s \exp(-\gamma \xi_v) \, dv \right| = 0 \qquad a.s.$$

and hence for every random variable $S$ with values in $[0, \infty]$,

$$(22) \quad \lim_{n \to \infty} \frac{K_n(S)}{n^\gamma \ell(n)} = \Gamma(1-\gamma) \int_0^S \exp(-\gamma \xi_v) \, dv \qquad a.s.$$

PROOF. The joint distribution of the two processes $(K_n(s), s \geq 0)$ and $(\xi_t, t \geq 0)$ is the same as if $V_1, V_2, \ldots$ were more specifically of the form $V_i = \inf\{v \geq 0 : e^{-\xi_v} < U_i\}$, where $U_1, U_2, \ldots$ is a sequence of independent uniform $(0,1)$ variables independent of $\xi$. Then $K_n(s)$ is, as in [26], the minimal number of open intervals of the form $(\exp(-\xi_v), \exp(-\xi_{v-}))$, $v \leq s$, containing $U_i, 1 \leq i \leq n$ [with $U_i \geq \exp(-\xi_s)$].

If $\mathbb{P}(S = s) = 1$ for some fixed $s \in [0, \infty]$, then the conclusion (22) is read from [26], Theorem 4.1, as indicated [26], Corollary 5.2, in the case $s = \infty$. The uniform convergence (21) follows by a standard pathwise argument, using the facts that the process $(K_n(s), s \geq 0)$ is increasing in $s$ for each $n$ and that the limit process $(\int_0^s \exp(-\gamma \xi_v) \, dv, s \geq 0)$ has continuous paths. □

PROOF OF PROPOSITION 7. Recall that the discrete tree $T_n$ is also considered as an $\mathbb{R}$-tree by letting its edge lengths all be 1, so we may consider reduced trees of the form $R(T_n, x_1, \ldots, x_k)$ where $x_1, \ldots, x_k$ are vertices of $T_n$. We let $R(T_n, B)$ be the reduced tree of $T_n$ spanned by the root and the vertices labeled by $B$. By exchangeability of the partition-valued process $(\Pi(t), t \geq 0)$, it is plain that $R(T_n, L_1^n, \ldots, L_k^n)$ has same law as $R(T_n, B)$ for any $B$ with $\#B = k$ and $n \geq \max B$. We are going to show that almost surely, for every finite $B \subset \mathbb{N}$,

$$(23) \qquad n^{-\gamma} \ell(n)^{-1} R(T_n, B) \xrightarrow[n \to \infty]{a.s.} \Gamma(1-\gamma) R(\mathcal{T}, \{L_i, i \in B\}).$$

Notice that the shape of $R(T_n, B)$ is exactly $T_B$, as in Definition 3, although the edge lengths are different from 1 in general.

Now, assume Hypothesis (H). Consider the case $k = 1$. By exchangeability, it is enough to discuss $B = \{1\}$. Let $D_1^n$ be the (combinatorial) distance between the root of $T_n$ and $\{1\}$. By construction of $T_n$, we see that $D_1^n - 1$ is the number of fragmentations that the block of $(\Pi^0|_n(t), t \geq 0)$ containing 1 undergoes from $[n]$ to $\{1\}$. Similarly,

$$(24) \quad D_1^n - 1 = \#\{L_1 \wedge L_i, 2 \leq i \leq n\} = \#\{d(\rho(\mathcal{T}), L_1 \wedge L_i), 2 \leq i \leq n\},$$

which is the number of branch points of $R(\mathcal{T}, L_1, \ldots, L_n)$ located on $[[\rho(\mathcal{T}), L_1]]$. Conditionally given $(\mathcal{T}, \mu)$ and $L_1$, the random variables $d(\rho(\mathcal{T}), L_1 \wedge L_i)$, $i \geq 2$, are independent and identically distributed with

$$\mathbb{P}(d(\rho(\mathcal{T}), L_1 \wedge L_i) > t | \mathcal{T}, \mu, L_1) = \mu(\mathcal{T}_t^{(1)}) = |\Pi_{(1)}(t)| = |\Pi^0_{(1)}(\eta_{(1)}(t))|,$$

where, according to (ii) in Proposition 6 and Lemma 5, the process $\eta_{(1)}$ is the inverse of the process

$$\eta_{(1)}^{-1} : t \longmapsto \int_0^t |\Pi^0_{(1)}(s)|^\gamma \, ds.$$

Moreover, $\xi := (-\log|\Pi^0_{(1)}(t)|, t \geq 0)$ is a pure-jump subordinator with Lévy measure defined in (18). Since the time-change $\eta_{(1)}$ is continuous and strictly increasing, we also see that

$$(25) \quad D_1^n - 1 = \#\{V_i, 2 \leq i \leq n\}, \qquad \text{where } V_i = \eta_{(1)}(d(\rho(\mathcal{T}), L_1 \wedge L_i))$$

so that, conditionally given $(\mathcal{T}, \mu)$ and $L_1$, the $V_i$ for $i \geq 2$ are independent and identically distributed with

$$\mathbb{P}(V_i > v | \mathcal{T}, \mu, L_1) = \exp(-\xi_v).$$

The desired conclusion that

$$n^{-\gamma} \ell(n)^{-1} D_1^n \xrightarrow[n \to \infty]{\text{a.s.}} \Gamma(1 - \gamma) \int_0^\infty \exp(-\gamma \xi_s) \, ds = \eta_{(1)}^{-1}(\infty) = D_1$$

is now read from (22) with $S = \infty$.

Next, assume that (23) holds for every $B$ with $\#B = k$. We show that it then holds for $\#B = k + 1$. Again by exchangeability, it is enough to discuss the case $B = [k + 1]$. Let $D_{[k+1]}$ be the first time $t$ when $[k + 1]$ is not included in a block of $\Pi(t)$ so that $D_{[k+1]}$ is, by definition, the length of the edge adjacent to the root in $R(\mathcal{T}, L_1, \ldots, L_{k+1})$, that is, the height of $L_1 \wedge L_2 \wedge \cdots \wedge L_{k+1}$. Similarly, we let $D^0_{[k+1]}$ be the analogous time, but for the process $(\Pi^0(t), t \geq 0)$. By the time-change correspondence between $\Pi$ and $\Pi^0$ (Lemma 5), if $\xi_t = -\log|\Pi^0_{(1)}(t)|, t \geq 0$, we know that

$$D_{[k+1]} = \int_0^{D^0_{[k+1]}} e^{-\gamma \xi_s} \, ds.$$

Let $D^n_{[k+1]}$ be the height of the first branch point in $R(T_n, \{1\}, \ldots, \{k+1\})$. Then $D^n_{[k+1]} - 1$ is the number of fragmentation events undergone by $\Pi|_n(t)$, $0 \leq t \leq D_{[k+1]}$. This is also the number of distinct branch points of $R(\mathcal{T}, L_1, \ldots, L_n)$ belonging to $R(\mathcal{T}, L_1, \ldots, L_{k+1})$ with height $\leq D_{[k+1]}$, that is,

$$D^n_{[k+1]} - 1 = \#\{L_1 \wedge L_i, k+2 \leq i \leq n, d(\rho(\mathcal{T}), L_1 \wedge L_i) \leq D_{[k+1]}\}$$

(one could use any $L_j$, $j \leq k+1$, instead of $L_1$ in this formula). By the same argument used for $k = 1$,

$$D^n_{[k+1]} - 1 = \#\{V_i, k+2 \leq i \leq n, V_i \leq D^0_{[k+1]}\}$$

$$\text{where } V_i = \eta_{(1)}(d(\rho(\mathcal{T}), L_1 \wedge L_i)),$$

as before. Formula (22) applied with $S = D^0_{[k+1]}$ now yields

$$(26) \quad n^{-\gamma}\ell(n)^{-1} D^n_{[k+1]} \xrightarrow[n\to\infty]{\text{a.s.}} \Gamma(1-\gamma) \int_0^{D^0_{[k+1]}} e^{-\gamma \xi_s}\, ds = \Gamma(1-\gamma) D_{[k+1]},$$

so the renormalized length of the root-edge of $R(T_n, \{1\}, \ldots, \{k+1\})$ converges to the length of the root-edge of $R(\mathcal{T}, L_1, \ldots, L_{k+1})$, up to the renormalization factor $\Gamma(1-\gamma)$.

Next, let $\pi = \Pi^0|_{k+1}(D^0_{[k+1]})$, with nonempty blocks $\pi_1, \ldots, \pi_r$. Recalling the notation of Sections 3.2 and 3.3.2, we have

$$R(T_n, [k+1]) = \langle R(T_n, \pi_1) - D^n_{[k+1]}, \ldots, R(T_n, \pi_r) - D^n_{[k+1]} \rangle_{D^n_{[k+1]}}$$

because $D^n_{[k+1]}$ is the height of the first branch point of $R(T_n, [k+1])$, while $\pi_i \subset [k+1]$. For the same reason,

$$R(\mathcal{T}, \{L_1, \ldots, L_{k+1}\})$$
$$= \langle R(\mathcal{T}, \{L_i, i \in \pi_1\}) - D_{[k+1]}, \ldots, R(\mathcal{T}, \{L_i, i \in \pi_r\}) - D_{[k+1]} \rangle_{D_{[k+1]}}.$$

Now, condition on the first split $\pi$. The conclusion follows from (26) and the induction hypothesis, which implies that for $1 \leq i \leq r$,

$$n^{-\gamma}\ell(n)^{-1} R(T_n, \pi_i) \xrightarrow[n\to\infty]{\text{a.s.}} \Gamma(1-\gamma) R(\mathcal{T}, \{L_j, j \in \pi_i\}).$$

This completes the proof under Hypothesis (H). Note, however, that the *joint* distribution of $(T_n)_{n \geq 1}$ as a sequence of $\Theta$-valued random variables is the same under Hypothesis (H) as in the apparently more general setting of Proposition 7. Since $\Theta$ is complete, we conclude that also in the setting of Proposition 7, there exists a tree $\mathcal{R}_k$ on the given probability space to which the rescaled $R(T_n, [k])$ converge a.s. □

4.2. *Tightness estimate.* The aim of this subsection is to prove the forthcoming tightness estimate (Proposition 9).

PROPOSITION 9. *For $k \leq n$, let*

$$\Delta(n,k) := \max_{1 \leq i \leq n} d_n(\{i\}, R(T_n, \{1\}, \ldots, \{k\})),$$

$d_n$ *being the metric associated with the tree $T_n$. Then, under the hypotheses of Theorem 2, for each $\eta > 0$,*

$$\lim_{k \to \infty} \limsup_{n \to \infty} \mathbb{P}\bigg(\frac{\Delta(n,k)}{\overline{\Lambda}(n^{-1})} > \eta\bigg) = 0.$$

Before we give the proof of this proposition, let us deduce Theorem 2.

PROOF OF THEOREM 2. First, assume Hypothesis (H). Fix $\varepsilon, \eta > 0$ and choose $k$ large enough that $\mathbb{P}(\Gamma(1-\gamma)d_{\mathrm{GH}}(\mathcal{R}_k, \mathcal{T}) > \eta) < \varepsilon$ (we know from [29] that $l_1$-representatives of $\mathcal{R}_k$ converge to $\mathcal{T}$ a.s. as $k \to \infty$; Hausdorff convergence in $l_1$ implies Gromov–Hausdorff convergence) and

$$\limsup_{n \to \infty} \mathbb{P}(d_{\mathrm{GH}}(R(T_n, \{1\}, \ldots, \{k\}), T_n) > \overline{\Lambda}(n^{-1})\eta) < \varepsilon$$

(such $k$ exists by Proposition 9). Then for $n$ sufficiently large,

$$\mathbb{P}(d_{\mathrm{GH}}(R(T_n, \{1\}, \ldots, \{k\}), T_n) > \overline{\Lambda}(n^{-1})\eta) < \varepsilon$$

and also

$$\mathbb{P}(d_{\mathrm{GH}}(R(T_n, \{1\}, \ldots, \{k\})/\overline{\Lambda}(n^{-1}), \Gamma(1-\gamma)\mathcal{R}_k) > \eta) < \varepsilon$$

since $R(T_n, \{1\}, \ldots, \{k\})/\overline{\Lambda}(n^{-1})$ converges a.s. to $\Gamma(1-\gamma)\mathcal{R}_k$ as $n \to \infty$ (see Proposition 7). Hence, for $n$ sufficiently large, $\mathbb{P}(d_{\mathrm{GH}}(T_n/\overline{\Lambda}(n^{-1}), \Gamma(1-\gamma)\mathcal{T}) > 3\eta) < 3\varepsilon$. This completes the proof for the setting of this section, where $T_n$, $n \geq 1$, are derived from an exchangeable sample of leaves $L_1, L_2, \ldots$ with directing measure $\mu$ of a given CRT $(\mathcal{T}, \mu)$. If we do not assume (H), then we argue, as at the end of the proof of Proposition 7, that for any probability space supporting $(T_n, n \geq 1)$, there exists a random $\mathbb{R}$-tree $\mathcal{T}_{(\gamma_\nu, \nu)}$ on the same probability space, to which the rescaled $T_n$ converge in probability. □

The proof of Proposition 9 which we postponed is given in Section 4.2.2, Section 4.2.1 being devoted to the proof of key intermediate results (Lemma 10 and its Corollary 11). We will work under Hypothesis (H), without loss of generality.

4.2.1. *A key lemma.* Throughout, we consider a fixed $\nu$. Implicitly, the constants appearing in this section may depend on $\nu$. Note that the conditions (6) and (7) satisfied by $\nu$ imply that the tail $\overline{\Lambda}$ of the Lévy measure $\Lambda$ [defined in (18)], that is, $\overline{\Lambda}(x) = \int_x^\infty \Lambda(dy)$, $x > 0$, satisfies both the regular variation condition $\overline{\Lambda}(x) \sim x^{-\gamma}\ell(1/x)$ as $x \to 0$ and $\int^\infty x^\rho \Lambda(dx) < \infty$. We may, and will, also assume, without loss of generality, that $0 < \rho < \gamma$. We claim that this implies the existence of some finite constant $C_\Lambda > 0$ such that

(27) $$\overline{\Lambda}(xy) \leq C_\Lambda \overline{\Lambda}(x) y^{-\rho} \qquad \text{for all } y \geq 1, 0 < x \leq 1.$$

To see this, choose $\delta < \gamma - \rho$ and note that Potter's theorem ([11], Theorem 1.5.6) implies the existence of some $X > 0$ such that for $x \in (0,1], y \geq 1$ with $xy \leq X$, we have

$$\overline{\Lambda}(xy)/\overline{\Lambda}(x) \leq 2y^{\delta-\gamma} \leq 2y^{-\rho}.$$

On the other hand, if $x \in (0,1]$ and $xy \geq X$, we have

$$\overline{\Lambda}(xy) \leq (xy)^{-\rho} \int_X^\infty z^\rho \Lambda(dz),$$

where the last integral is finite, while $x^{-\rho} \leq C'_\Lambda \overline{\Lambda}(x)$ for some constant $C'_\Lambda > 0$ because $x^{-\rho}/\overline{\Lambda}(x)$ is regularly varying with exponent $\gamma - \rho > 0$ at 0 and is hence bounded on $(0,1]$. The estimate (27) will be useful in the sequel.

Let $H_n$ be the height of the tree $T_n$, that is, $H_n := \max_{1 \leq i \leq n} D_i^n$, where $D_i^n$ denotes the height of the leaf $\{i\}$ (i.e., its distance to the root) in the tree $T_n$.

LEMMA 10. *There exists a random variable $X_\infty$, with positive moments of all orders, such that, for all $p \geq 2/\gamma$, there exists a constant $C_p$ such that, for all $x \geq 1$ and all integers $n$,*

$$\mathbb{P}(H_n > (1+x)2X_\infty \overline{\Lambda}(n^{-1})) \leq \frac{C_p}{x^p}.$$

COROLLARY 11. *For all $a > 0$ and $p \geq 2/\gamma$, there exists some constant $C_{p,a}$ such that, for all $x \geq 1$ and all integers $n$,*

$$\mathbb{P}(H_n > ax\overline{\Lambda}(n^{-1})) \leq \frac{C_{p,a}}{x^p}.$$

PROOF. We simply use the fact that

$$\mathbb{P}(H_n > ax\overline{\Lambda}(n^{-1})) \leq \mathbb{P}(H_n > ax\overline{\Lambda}(n^{-1}), ax \geq (1+\sqrt{x})2X_\infty)$$
$$+ \mathbb{P}((1+\sqrt{x})2X_\infty > ax),$$

then bound the right-hand side from above side using the upper bound of Lemma 10 for the first probability and the fact that $E[X_\infty^{2p}]$ is finite for the second probability. □

The main idea needed to prove Lemma 10 is to transfer the problem on the tail of $H_n$ onto a problem on the tail of $D_1^n$, using $H_n = \max_{1 \le i \le n} D_i^n$. Indeed, for all (random) sequences $(X_i)_{i \ge 1}$ such that the random variables $(D_i^n, X_i)$, $1 \le i \le n$, are identically distributed, one has

$$\mathbb{P}(H_n > X_\infty x) \le n \mathbb{P}(D_1^n > X_1 x) \qquad \forall x \ge 0,$$

where $X_\infty := \sup_{i \ge 1} X_i$. Therefore, it is sufficient to find random variables $X_i$, $i \ge 1$, whose supremum possesses moments of all positive orders and then a suitable upper bound for the tail of $D_1^n$ to conclude. This is the goal of the remainder of this subsection. Define $X_i$ by

$$X_i := (1 + A_\gamma) C_\Lambda \sum_{k=0}^\infty \exp(-\rho \xi_k^i) + 1,$$

where

$$A_\gamma := 2 \sum_{k=1}^\infty \frac{(k+1)^{\sqrt{\gamma}}}{k(k+1)} < \infty,$$

since $\gamma \in (0,1)$, and $\xi^i$ is the subordinator describing the evolution of the sizes of the blocks containing $i$ in the fragmentation $\Pi^0$, as explained in Lemma 5. Clearly, $(D_i^n, X_i)$, $1 \le i \le n$, are identically distributed (by exchangeability) and

$$X_\infty = \sup_{i \ge 1} X_i \le (1 + A_\gamma) C_\Lambda (1 + \zeta_\rho) + 1,$$

where

$$\zeta_\rho := \sup_{i \ge 1} \int_0^\infty \exp(-\rho \xi_t^i) \, dt$$

is (in distribution) the first time at which a self-similar fragmentation with parameters $(-\rho, \nu)$ reaches the trivial partition $\{\{1\}, \{2\}, \ldots\}$ (in others words, it is the height of the associated fragmentation tree). It was proven in [28] (Proposition 14) that $\zeta_\rho$ (hence $X_\infty$) has exponential moments. Lemma 10 is therefore an immediate consequence of the following result.

LEMMA 12. *For all $p \ge 0$, there exists a constant $C_p'$ such that for all $x \ge 1$ and all integers $n$,*

$$\mathbb{P}(D_1^n > (1+x) 2 X_1 \overline{\Lambda}(n^{-1})) \le \frac{C_p'}{x^p n^{\gamma p - 1}}.$$

The remainder of this subsection is devoted to the proof of this lemma. To simplify the notation, we omit the index 1 wherever we can (i.e., $\xi$ now stands for $\xi^1$, $X$ for $X_1$). We also set $D_n := D_1^{n+1} - 1$ for the number of internal vertices between the root and leaf $\{1\}$ in $T_{n+1}$. Since $D_1^{n+1} \leq 2D_n$ and $\overline{\Lambda}(n^{-1}) \leq \overline{\Lambda}((n+1)^{-1})$, the upper bound stated in Lemma 12 is a consequence of the existence of some constant $C_p''$ such that for all $x \geq 1$ and all integers $n$,

$$(28) \qquad \mathbb{P}(D_n > (1+x)X\overline{\Lambda}(n^{-1})) \leq \frac{C_p''}{x^p n^{\gamma p - 1}}.$$

To prove this latter inequality, we proceed in three steps.

Let $N_x(s,t)$ denote the number of jumps of $\xi$ of size at least $x$ in the time interval $[s,t]$, $\tilde{N}_x(s,t)$ denote the number of jumps of $1 - \exp(-\xi)$ of size at least $x$ in the same time interval and $\tilde{N}_x := \tilde{N}_x(0, \infty)$.

**Step 1. Large deviations for $\tilde{N}_x$.** The regular variation of $\overline{\Lambda}$ at 0 ensures that $\tilde{N}_x \sim \overline{\Lambda}(x)D$ a.s. as $x \to 0$, where $D = \int_0^\infty \exp(-\gamma \xi_t)\,dt$ (Theorem 5.1, [26]). The goal of this first step is to give some kind of large deviations result on this convergence.

LEMMA 13.  *For all $x > 0$ and $0 < y \leq 1$,*

$$\mathbb{P}\left(\tilde{N}_y > (1+x)C_\Lambda \sum_{i=0}^\infty (\exp(-\rho \xi_i))\overline{\Lambda}(y)\right) \leq \exp(-a_x \overline{\Lambda}(y)),$$

*where $a_x := (1+x)\ln(1+x) - x > 0$.*

PROOF. Let $\mathcal{F}_t$ denote the $\sigma$-field generated by $\xi$ until time $t$ and $\mathcal{F}$ the one generated by $\xi$, and observe that

$$\tilde{N}_y = \sum_{i=0}^\infty \tilde{N}_y(i, i+1) \leq \sum_{i=0}^\infty N_{y\exp(\xi_i)}(i, i+1).$$

Conditional on $\mathcal{F}_i$, $N_{y\exp(\xi_i)}(i, i+1)$ is a Poisson random variable with mean $\overline{\Lambda}(y\exp(\xi_i))$. But for any Poisson random variables $P$ with mean $\lambda$, one has

$$\mathbb{E}[\exp(tP - (1+x)t\lambda)] = \exp((\exp(t) - 1 - (1+x)t)\lambda) \qquad \forall t \in \mathbb{R}.$$

In particular, when $t = \ln(1+x)$, $\exp(t) - 1 - (1+x)t = -a_x < 0$ and the expectation is smaller than 1. Hence, for all $n \in \mathbb{N}$, using (27) for the first inequality, we get, for all $y \leq 1$,

$$\mathbb{P}\left(\sum_{i=0}^n N_{y\exp(\xi_i)}(i, i+1) \geq (1+x)C_\Lambda \sum_{i=0}^n \exp(-\rho \xi_i)\overline{\Lambda}(y)\right)$$

$$\leq \mathbb{P}\left(\sum_{i=0}^n N_{y\exp(\xi_i)}(i, i+1) \geq (1+x)\sum_{i=0}^n \overline{\Lambda}(y\exp(\xi_i))\right)$$

$$\leq \mathbb{E}\left[\exp\left(t\left(\sum_{i=0}^{n}(N_{y\exp(\xi_i)}(i,i+1)-(1+x)\overline{\Lambda}(y\exp(\xi_i)))\right)\right)\right]$$

$$\leq \mathbb{E}\left[\exp\left(t\left(\sum_{i=0}^{n-1}\cdots\right)\right)\mathbb{E}[\exp(t(N_{y\exp(\xi_n)}(n,n+1)$$

$$-(1+x)\overline{\Lambda}(y\exp(\xi_n))))|\mathcal{F}_n]\right]$$

$$\leq \cdots \leq \exp(-a_x\overline{\Lambda}(y)),$$

the last line being obtained by induction: at each step but the last, we use the upper bound 1 for the (conditional) expectation and for the last step, we use the upper bound $\exp(-a_x\overline{\Lambda}(y))$ for the expectation $\mathbb{E}[\exp(t(N_y(0,1) - (1+x)\overline{\Lambda}(y)))]$. It remains to let $n \to \infty$ in the first probability involved in the above sequence of inequalities and to use Fatou's lemma. □

**Step 2. Large deviations for $\mathbb{E}[D_n|\mathcal{F}]$.** We now establish a result similar to the required inequality (28), but for the quantity $\mathbb{E}[D_n|\mathcal{F}]$, where $\mathcal{F} = \mathcal{F}_\infty$ is the $\sigma$-field generated by the whole subordinator $\xi$ [recall that we work under Hypothesis (H)].

LEMMA 14. *Let $B_\gamma := \sum_{k=1}^\infty \exp(-4^{-1}a_1 k^{\gamma/2})$ with $a_1 = 2\ln 2 - 1$. Then for all $x \geq 1$ and all integers $n$ large enough,*

$$\mathbb{P}(\mathbb{E}[D_n|\mathcal{F}] > (1+x)(X-1)\overline{\Lambda}(n^{-1})) \leq (1+B_\gamma)\exp(-4^{-1}a_1 x\overline{\Lambda}(n^{-1})).$$

PROOF. According to the formula (4) of [26],

$$\mathbb{E}[D_n|\mathcal{F}] = n\int_0^1 (1-y)^{n-1}\tilde{N}_y\,dy \leq \tilde{N}_{1/n} + n\int_0^{1/n}\tilde{N}_y\,dy.$$

Hence, setting $S := C_\Lambda \sum_{i=0}^\infty \exp(-\rho\xi_i)$,

$$\mathbb{P}(\mathbb{E}[D_n|\mathcal{F}] > (1+x)(1+A_\gamma)S\overline{\Lambda}(n^{-1}))$$
$$\leq \mathbb{P}(\tilde{N}_{1/n} > (1+x)S\overline{\Lambda}(n^{-1}))$$
$$+ \mathbb{P}\left(n\int_0^{1/n}\tilde{N}_y\,dy > (1+x)A_\gamma S\overline{\Lambda}(n^{-1})\right).$$

The first probability in the right-hand side is smaller than $\exp(-a_x\overline{\Lambda}(n^{-1}))$, according to Lemma 13. To bound the second probability, we use $n\int_{1/(k+1)n}^{1/kn}\tilde{N}_y\,dy \leq \tilde{N}_{1/(n(k+1))}\frac{1}{k(k+1)}$, which gives

$$\mathbb{P}\left(n\int_0^{1/n}\tilde{N}_y\,dy > A_\gamma(1+x)S\overline{\Lambda}(n^{-1})\right)$$

$$\leq \sum_{k=1}^{\infty} \mathbb{P}(\tilde{N}_{1/n(k+1)} > 2(k+1)^{\sqrt{\gamma}}(1+x)S\overline{\Lambda}(n^{-1})).$$

Since $\overline{\Lambda}$ is regularly varying at 0 with index $-\gamma$, we have, provided $n$ is large enough, that $\overline{\Lambda}(n^{-1})(k+1)^{\gamma/2} \leq 2\overline{\Lambda}(((k+1)n)^{-1})$ and $\overline{\Lambda}(((k+1)n)^{-1}) \leq 2\overline{\Lambda}(n^{-1})(k+1)^{\sqrt{\gamma}}$ for all $k \geq 1$ (to see this, use, e.g., Potter's theorem, Theorem 1.5.6, [11]). Combined with Lemma 13, this implies that the above sum of probabilities is smaller than

$$\sum_{k=1}^{\infty} \exp(-a_x\overline{\Lambda}(((k+1)n)^{-1})) \leq \sum_{k=1}^{\infty} \exp(-2^{-1}a_x\overline{\Lambda}(n^{-1})(k+1)^{\gamma/2}).$$

Last, the exponential in the latter sum can be split in two, using $(k+1)^{\gamma/2} \geq 2^{-1}(k^{\gamma/2}+1)$, to get the upper bound

$$\exp(-4^{-1}a_x\overline{\Lambda}(n^{-1})) \sum_{k=1}^{\infty} \exp(-a_x 4^{-1}\overline{\Lambda}(n^{-1})k^{\gamma/2}),$$

which is smaller than $\exp(-4^{-1}a_1 x\overline{\Lambda}(n^{-1}))B_\gamma$ for all $x \geq 1$ ($a_x \geq a_1 x$ for $x \geq 1$) and $n$ large enough. $\square$

**Step 3. Proof of inequality (28).** To start with, fix $x \geq 1$, $n \in \mathbb{N}$, and note that

$$\begin{aligned}(29)\quad \mathbb{P}(D_n > (1+x)X\overline{\Lambda}(n^{-1})) &\leq \mathbb{P}(\mathbb{E}[D_n|\mathcal{F}] > (1+x)(X-1)\overline{\Lambda}(n^{-1})) \\ &\quad + \mathbb{P}(D_n - \mathbb{E}[D_n|\mathcal{F}] > (1+x)\overline{\Lambda}(n^{-1})).\end{aligned}$$

Lemma 14 gives an upper bound for the first probability, provided $n$ is large enough. To get an upper bound for the second probability, we use a result on urn models (Devroye [12], Section 6) which ensures that

$$\mathbb{P}(D_n - \mathbb{E}[D_n|\mathcal{F}] > y|\mathcal{F}) \leq \exp\left(-\frac{y^2}{2\mathbb{E}[D_n|\mathcal{F}] + 2y/3}\right) \qquad \forall y \geq 0,\ n \in \mathbb{N}.$$

This implies that for all $m \geq 1$, there exists some deterministic constant $B_m$ depending only on $m$ such that

$$\begin{aligned}\mathbb{P}(D_n &- \mathbb{E}[D_n|\mathcal{F}] > (1+x)\overline{\Lambda}(n^{-1})|\mathcal{F}) \\ &\leq B_m\left(\frac{\mathbb{E}[D_n|\mathcal{F}] + (1+x)\overline{\Lambda}(n^{-1})}{((1+x)\overline{\Lambda}(n^{-1}))^2}\right)^m \\ &\leq 2^{m-1}B_m \frac{(\mathbb{E}[D_n|\mathcal{F}])^m + ((1+x)\overline{\Lambda}(n^{-1}))^m}{((1+x)\overline{\Lambda}(n^{-1}))^{2m}} \\ &\leq 2^{m-1}B_m \frac{\mathbb{E}[D_n^m|\mathcal{F}] + ((1+x)\overline{\Lambda}(n^{-1}))^m}{((1+x)\overline{\Lambda}(n^{-1}))^{2m}},\end{aligned}$$

the last line being obtained by Jensen's inequality. We then take expectations on both sides of the resulting inequality. Theorem 6.3 of [26] ensures that $\mathbb{E}[D_n^m] \sim (\overline{\Lambda}(n^{-1}))^m$ (up to a constant). Therefore, we have

$$(30) \quad \mathbb{P}(D_n - \mathbb{E}[D_n|\mathcal{F}] > (1+x)\overline{\Lambda}(n^{-1})) \leq B_{m,\Lambda}((1+x)\overline{\Lambda}(n^{-1}))^{-m},$$

where $B_{m,\Lambda}$ depends only on $m$ and $\Lambda$.

Next, recall the upper bound given by Lemma 14 for the first probability involved in the right-hand side of (29). Together with the upper bound (30), it leads to the existence of $B'_{m,\Lambda}$ such that

$$\mathbb{P}(D_n > (1+x)X\overline{\Lambda}(n^{-1})) \leq B'_{m,\Lambda} x^{-m} (\overline{\Lambda}(n^{-1}))^{-m}$$

for all $x \geq 1$ and $n$ large enough, say $n \geq n_0$. Since $\overline{\Lambda}(n^{-1}) \sim n^\gamma \ell(n)$ when $n \to \infty$, this upper bound is, in turn, bounded from above by $x^{-m} n^{1-\gamma m}$, up to some constant, which is the required result (28).

Finally, inequality (28) is also true when $n \leq n_0$ (for all $x \geq 1$) since $D_n \leq n \leq n_0$ and $X \geq 1$, and therefore the probability $\mathbb{P}(D_n > (1+x)X\overline{\Lambda}(n^{-1}))$ is null whenever $1 + x \geq n_0(\overline{\Lambda}(1))^{-1}$.

4.2.2. *Proof of Proposition 9.* The crucial point is that

$$\Delta(n,k) = \max_{j \geq 1} \overline{H}_{n_j^{k,n}},$$

where the $n_j^{k,n}$ and $\overline{H}_{n_j^{k,n}}$, $1 \leq k \leq n$, $j \geq 1$, are defined as follows. Let $\Pi_{(i)}(t)$ denote the block of $\Pi(t)$ containing $i$, $i \geq 1$. Then for all $k \geq 1$, introduce

$$t_i^k := \inf\{t \geq 0 : \Pi_{(i)}(t) \cap [k] = \varnothing\},$$

the first time at which the fragment containing $i$ is disjoint from $[k]$ (in particular, $t_i^k = \infty$ for $1 \leq i \leq k$). For all $t \geq 0$, the collection of blocks $(\Pi_{(i)}(t_i^k + t), i \geq k+1)$ induces a partition, denoted $\Pi(t^k + t)$, of $\mathbb{N} \setminus [k]$ and each $\Pi_j(t^k + t)$ admits asymptotic frequencies, as $\Pi(t^k + t)$ is an exchangeable partition of $\mathbb{N}\setminus[k]$. We call $n_j^{k,n}$ the cardinality of $\Pi_j(t^k) \cap [n]$ and $\lambda_j^k$ the a.s. limit of $n_j^{k,n}/n$ as $n \to \infty$. Clearly, $\lambda_{\max}^k := \max_{j \geq 1} \lambda_j^k \to 0$ a.s. as $k \to \infty$.

Then let $\mathcal{G}(k)$ be the $\sigma$-field generated by $\Pi(t^k)$. In the terminology of Bertoin ([10], Definition 3.4), the sequence $(t_i^k, i \in \mathbb{N})$ is a *stopping line* and, as such, satisfies the *extended branching property* ([10], Lemma 3.14) which ensures that given $\mathcal{G}(k)$, the process $(\Pi(t^k + t), t \geq 0)$ is a fragmentation process starting from $\Pi(t^k)$. This implies that given $\mathcal{G}(k)$, the discrete fragmentation trees, with respectively $n_1^{k,n}, n_2^{k,n}, \ldots$ leaves, associated with the fragmentations of the blocks $\Pi_j(t^k)$, $j \geq 1$, evolve independently as $n \to \infty$ with laws respectively distributed as $T_{n_j^{k,n}}$, $j \geq 1$. In particular, given $\mathcal{G}(k)$, the respective heights of those trees, $\overline{H}_{n_j^{k,n}}$, $j \geq 1$, are independent and distributed as $H_{n_j^{k,n}}$, $j \geq 1$.

Let $\eta > 0$. We now turn back to our goal, which is to prove that
$$\lim_{k \to \infty} \liminf_{n \to \infty} \mathbb{P}(\Delta(n,k) \leq \eta \overline{\Lambda}(n^{-1})) = 1.$$

Note that first applying dominated convergence (for the limit in $k$, everything is bounded by 1) and then Fatou's lemma (for the lim inf in $n$), it is sufficient to show that
$$\lim_{k \to \infty} \liminf_{n \to \infty} \mathbb{P}(\Delta(n,k) \leq \eta \overline{\Lambda}(n^{-1}) | \mathcal{G}(k)) \to 1 \quad \text{a.s.}$$

According to the discussion above,
$$\mathbb{P}(\Delta(n,k) \leq \eta \overline{\Lambda}(n^{-1}) | \mathcal{G}(k)) = \prod_{j \geq 1} \mathbb{P}(\overline{H}_{n_j^{k,n}} \leq \eta \overline{\Lambda}(n^{-1}) | \mathcal{G}(k))$$

and our goal turns into the proof of
$$\lim_{k \to \infty} \liminf_{n \to \infty} \sum_{j \geq 1} \ln(1 - \mathbb{P}(\overline{H}_{n_j^{k,n}} > \eta \overline{\Lambda}(n^{-1}) | \mathcal{G}(k))) = 0.$$

For the rest of the argument, we may consider that $n_j^{k,n}$, $\lambda_j^k$, $j \geq 1$, are deterministic and drop the conditioning on $\mathcal{G}(k)$ from the notation. Let $p > \max(\rho^{-1}, 2/\gamma)$. By inequality (27), for all $j, k, n \geq 1$ such that $n_j^{k,n} \neq 0$,
$$C_\Lambda \overline{\Lambda}(n^{-1}) \geq \left(\frac{n}{n_j^{k,n}}\right)^\rho \overline{\Lambda}((n_j^{k,n})^{-1}).$$

Corollary 11 then ensures that
$$\mathbb{P}(H_{n_j^{k,n}} > \eta \overline{\Lambda}(n^{-1})) \leq C_{p,\Lambda,\eta} \left(\frac{n_j^{k,n}}{n}\right)^{p\rho},$$
where $C_{p,\Lambda,\eta}$, depends only on $p$, $\Lambda$ and $\eta$, for all $i, k, n \geq 1$, with the convention $H_0 := 0$.

In the rest of the proof, we choose $k$ large enough, say $k \geq k_0$, so that $\lambda_{\max}^k \leq (2(2C_{p,\Lambda,\eta})^{1/p\rho})^{-1}$. Then consider some integer $j_k$ such that $\sum_{j \geq j_k} \lambda_j^k \leq \lambda_{\max}^k$. Since $n_j^{k,n}/n \to \lambda_j^k$ as $n \to \infty$ for all $j \geq 1$ and also $\sum_{j \geq j_k} n_j^{k,n}/n \to \sum_{j \geq j_k} \lambda_j^k$, there exists an integer $n_k$ such that for all $n \geq n_k$, $n_j^{k,n}/n \leq 2\lambda_j^k$, $1 \leq j < j_k$, and $\sum_{j \geq j_k} n_j^{k,n}/n \leq 2\lambda_{\max}^k$. In particular, $n_j^{k,n}/n \leq 2\lambda_{\max}^k$ for all $j \geq 1$. Consequently, using the fact that $|\ln(1-x)| \leq 2x$ when $0 < x \leq 1/2$, we have for all $n \geq n_k$,
$$\sum_{j \geq 1} |\ln(1 - \mathbb{P}(H_{n_j^{k,n}} > \eta \overline{\Lambda}(n^{-1})))| \leq 2C_{p,\Lambda,\eta} \sum_{j \geq 1} \frac{(n_j^{k,n})^{p\rho}}{n^{p\rho}}$$

$$\leq 2C_{p,\Lambda,\eta}\left(\sum_{j=1}^{j_k}\frac{(n_j^{k,n})^{p\rho}}{n^{p\rho}}+(2\lambda_{\max}^k)^{p\rho}\right).$$

The parenthesis in the upper bound converges to $(\sum_{i=1}^{i_k}(\lambda_i^k)^{p\rho}+(2\lambda_{\max}^k)^{p\rho})$ as $n\to\infty$, which is smaller than $(\lambda_{\max}^k)^{p\rho-1}(1+2^{p\rho})$ (since $\sum_{i=1}^{i_k}\lambda_i^k\leq 1$). The result follows since $\lambda_{\max}^k\to 0$ as $k\to\infty$.

4.3. *Height functions.* The aim of this section is to provide an analog of Theorem 2 for a family of functions coding the heights of leaves in ordered versions of the trees. In the special case of beta-splitting models, this convergence of leaf-height functions was suggested, but not proven, by Aldous [1].

The ordered version of $T_n$ is obtained by putting the set of children of every nonleaf vertex of $T_n$ in exchangeable random order, independently over distinct vertices and conditionally on $T_n$. This is usually achieved by taking (rooted) planar embeddings of the trees, where the order among children of a vertex is read from the clockwise ordering of edges going from the vertex to its children. We then define the order $\preceq_n$ as a linear order on the leaves $\{1\},\ldots,\{n\}$ of $T_n$ by saying that $\{i\}\preceq_n\{j\}$ if the subtree pending from the most recent common ancestor $\{i\}\wedge\{j\}$ of $\{i\}$ and $\{j\}$ that contains $\{i\}$ comes before the subtree pending from $\{i\}\wedge\{j\}$ that contains $\{j\}$.

If $(T_n, n\geq 1)$ is a strongly consistent family of trees, we also want the orders $(\preceq_n, n\geq 1)$ to satisfy a consistency property, namely, the restriction of $\preceq_{n+1}$ to $\{1\},\{2\},\ldots,\{n\}$ is $\preceq_n$. With our interpretation of ordered trees as planar embeddings, this means that the embeddings are drawn consistently. This can be achieved inductively as follows, starting from $\preceq_1$, the trivial order on $\{\{1\}\}$. Suppose we are given $T_{n+1}$ and $\preceq_n$. Denote by $b(\{n+1\})$ the father of $\{n+1\}$ in $T_{n+1}$. For any nonleaf vertex $v$ of $T_{n+1}$ distinct from $b(\{n+1\})$, the children of $v$ are ordered in the same way for $T_n, \preceq_n$. Hence, the restriction to $\{1\},\ldots,\{n\}$ of $\preceq_{n+1}$ must be $\preceq_n$.

Next, two possibilities occur: either $b(\{n+1\})$ was already a vertex of $T_n$ or $b(\{n+1\})$ is a newly added vertex in $T_{n+1}$ with two offspring.

- If $b(\{n+1\})$ is a vertex of $T_n$ with $r$ children ordered as $c_1,\ldots,c_r$, we let $\{n+1\}$ be the $j$th son of $b(\{n+1\})$ in $T_{n+1}$, $1\leq j\leq r+1$, with equal probability $1/(r+1)$, and the order of the other children is preserved.
- Otherwise, $b(\{n+1\})$ must have a unique son $c$ besides $\{n+1\}$ in $T_{n+1}$ and we let $\{n+1\}$ be placed before or after $c$ with equal probability $1/2$.

Note that $\preceq_n$ naturally extends to a linear order on $T_n$ by letting $v\preceq_n w$ if either $v$ is an ancestor of $w$ or $v\wedge w=\{i\}\wedge\{j\}$ for some leaves $\{i\},\{j\}$ such that $\{i\}\preceq_n\{j\}$. This corresponds to the usual depth-first search order for rooted planar trees.

For each $n \geq 1$, we associate with the ordered tree $T_n^{\mathrm{ord}} = (T_n, \preceq_n)$ its *leaf-height function* $h_n$, defined on $[0,1]$ by $h_n(0) := 0$, $h_n(1) := 0$ and, for $1 \leq i \leq n$,

$$h_n\left(\frac{i}{n+1}\right) := \text{height of the } i\text{th leaf (in the left-to-right ordering)},$$

with linear interpolation. In general, the leaf-height function does not encode the full shape of the discrete tree and, more precisely, leaves some ambiguity where there are multiple branch points (e.g., the two possible unlabeled ordered trees with five leaves, all at distance 3 from the root vertex, are not distinguished by the leaf-height process).

Similarly, but fully encoding, a continuous *height* function $h:[0,1] \to \mathbb{R}^+$, $h(0) = h(1) = 0$, can be associated with the limiting fragmentation tree $\mathcal{T}$. Roughly, the construction of $h$ proceeds as follows (for details, we refer to Theorem 3 and Section 4.1 of [29], where it is more precisely proved that any fragmentation tree with an *infinite* dislocation measure—which is the case here—can be encoded into such continuous function). For each $k, n$ such that $k \leq n$, let $I_k^n \in \{1, \ldots, n\}$ be the position of the leaf $\{k\}$ among the leaves of $T_n$, with respect to the left-to-right ordering $\preceq_n$. Then define

$$U_k := \lim_{n \to \infty} \frac{I_k^n}{n+1}.$$

These limits exist a.s. and the $U_k$, $k \geq 1$, are i.i.d. uniformly distributed on $[0,1]$. The height function $h$ is then defined on $\{U_k, k \geq 1\}$ by $h(U_k) := $ height of $\{k\}$ in $\mathcal{T}$ and its definition can be extended continuously to $[0,1]$. The tree $\mathcal{T}$ can be recovered from $h$: it is isometric to the quotient space $([0,1], \overline{d})/\sim$, where $\overline{d}(x,y) := h(x) + h(y) - 2\inf_{z \in [x,y]} h(z)$ and $x \sim y \Leftrightarrow \overline{d}(x,y) = 0$. An order $\preceq$ on the leaves of $\mathcal{T}$ is then implicitly given by the natural order on $[0,1]$: let $x, y \in [0,1]$; if their images $\overline{x}, \overline{y}$ by projection on the quotient space are leaves, then $x \leq y \Leftrightarrow \overline{x} \preceq \overline{y}$. Further, according to Theorem 4 of [29], the function $h$ is a.s. Hölder-continuous of any order $\theta < \gamma$, but not of order $\theta > \gamma$ when $\nu$ integrates $s_1^{-1}$.

The a.s. convergence in Proposition 7 gives us a first connection between $h_n$ and $h$; namely, for all $k$,

$$\begin{aligned}
& n^{-\gamma} \ell(n)^{-1} \left( h_n\left(\frac{I_1^n}{n+1}\right), \ldots, h_n\left(\frac{I_k^n}{n+1}\right) \right) \\
& \qquad \xrightarrow[n \to \infty]{\mathrm{a.s.}} \Gamma(1-\gamma)(h(U_1), \ldots, h(U_k)).
\end{aligned}$$
(31)

More precisely, the following holds.

THEOREM 15. *In the situation of Theorem 2,*

$$\left(\frac{h_n(t)}{n^\gamma \ell(n)}\right)_{0 \leq t \leq 1} \xrightarrow[n \to \infty]{(p)} (\Gamma(1-\gamma) h(t))_{0 \leq t \leq 1}$$

*for the uniform norm on the space of continuous functions on $[0,1]$.*

PROOF. Let $\overline{h}_n := h_n/n^\gamma \ell(n)\Gamma(1-\gamma)$ and note that the convergences (31) imply that the only possible uniform limit (in distribution) for subsequences of $\overline{h}_n$ is $h$. Similarly to the proof of Theorem 2, we can strengthen (31) into convergence in probability for the uniform norm, by using a certain uniform estimate. This is inspired by a tightness estimate used the proof of Aldous [4], Theorem 20 for convergence of contour functions.

Fix $k \leq n$ and consider the order statistics $I_{(i)}^{n,k}, 1 \leq i \leq k$, of $I_k^n, 1 \leq i \leq k$. Also let $I_{(0)}^{n,k} := 0$, $I_{(k+1)}^{n,k} := n+1$. Then introduce

$$w_k^0(\overline{h}_n) := \max_{0 \leq i \leq k} \sup_{t \in [I_{(i)}^{n,k}/(n+1), I_{(i+1)}^{n,k}/(n+1)]} |\overline{h}_n(t) - \overline{h}_n(I_{(i)}^{n,k}/(n+1))|.$$

Our goal is to prove that

$$(32) \qquad \lim_{k \to \infty} \limsup_{n \to \infty} \mathbb{P}(w_k^0(\overline{h}_n) > \eta) = 0 \qquad \forall \eta > 0,$$

which is the analog of formula (30) of Aldous [4], Theorem 20, with $\alpha = 0$ there. Following the last lines of the proof of Aldous, one sees that (32) implies the tightness of $(\overline{h}_n, n \geq 1)$.

To get (32), first note that

$$w_k^0(\overline{h}_n) \leq \max_{0 \leq i \leq k} \left| \max_{t \in [I_{(i)}^{n,k}/(n+1), I_{(i+1)}^{n,k}/(n+1)]} \overline{h}_n(t) - \min_{t \in [I_{(i)}^{n,k}/(n+1), I_{(i+1)}^{n,k}/(n+1)]} \overline{h}_n(t) \right|$$

$$\leq \max_{0 \leq i \leq k} |\overline{d}_n(\{q_{\max}^{n,k,i}\}, \{q_{\min}^{n,k,i}\})|,$$

where $\overline{d}_n$ is the metric associated with the real trees $\mathcal{T}_n/(n^\gamma \ell(n)\Gamma(1-\gamma))$ and $\{q_{\max}^{n,k,i}\}$ is the leaf of $\mathcal{T}_n$ that has the highest height among the leaves $\{j\}$ of $\mathcal{T}_n$ such that $\{I_{(i)}^{n,k}\} \preceq_n \{j\} \preceq_n \{I_{(i+1)}^{n,k}\}$, where, by convention, both $\{I_{(0)}^{n,k}\}$ and $\{I_{(k+1)}^{n,k}\}$ denote the root $\rho(\mathcal{T})$ of $\mathcal{T}$. Similarly, among these leaves, $\{q_{\min}^{n,k,i}\}$ is the one that has the lowest height. Then define $v_{\max}^{n,k,i}$ in $\mathcal{R}_k^n := \mathcal{R}(\mathcal{T}_n, \{1\}, \ldots, \{k\})$ by

$$\overline{d}_n(\{q_{\max}^{n,k,i}\}, v_{\max}^{n,k,i}) = \overline{d}_n(\{q_{\max}^{n,k,i}\}, \mathcal{R}_k^n)$$

and define similarly $v_{\min}^{n,k,i}$. Now, fix $\varepsilon, \eta > 0$. Proposition 9 ensures that for $k$ large enough and then for $n$ sufficiently large,

$$\mathbb{P}\left( \max_{0 \leq i \leq k} \overline{d}_n(\{q_{\max}^{n,k,i}\}, v_{\max}^{n,k,i}) > \eta \right) \leq \varepsilon$$

and

$$\mathbb{P}\left( \max_{0 \leq i \leq k} \overline{d}_n(\{q_{\min}^{n,k,i}\}, v_{\min}^{n,k,i}) > \eta \right) \leq \varepsilon.$$

On the other hand, $\overline{d}_n(v_{\max}^{n,k,i}, v_{\min}^{n,k,i}) \leq \overline{d}_n(\{I_{(i)}^{n,k}\}, \{I_{(i+1)}^{n,k}\})$ and, using Proposition 7,

$$\max_{0 \leq i \leq k} \overline{d}_n(\{I_{(i)}^{n,k}\}, \{I_{(i+1)}^{n,k}\}) \to \max_{0 \leq i \leq k} d(L_{(i)}^k, L_{(i+1)}^k) \qquad \text{a.s. as } n \to \infty,$$

where $L_{(1)}^k \preceq \cdots \preceq L_{(k)}^k$ denotes the $\preceq$-ordered sequence of leaves $\{1\}, \ldots, \{k\}$ in $\mathcal{T}$ and $L_{(0)}^k := L_{(k+1)}^k := \rho(\mathcal{T})$. Finally, it is not hard to check that the compactness of $\mathcal{T}$ and its *ordered leaf-density* (informally, this means that the leaves are dense with respect to the order $\preceq$; see [29], Section 4.1, for precise details) imply that $\max_{0 \leq i \leq k} d(L_{(i)}^k, L_{(i+1)}^k) \to 0$ a.s. as $k \to \infty$. Therefore, for $k$ large enough and then for $n$ sufficiently large,

$$\mathbb{P}(w_k^0(\overline{h}_n) > 3\eta) \leq 3\varepsilon,$$

hence (32).

With this available, we just write

$$\sup_{0 \leq t \leq 1} |\overline{h}_n(t) - h(t)|$$

$$\leq w_k^0(\overline{h}_n) + \max_{0 \leq i \leq k} \sup_{t \in [I_{(i)}^{n,k}/(n+1), I_{(i+1)}^{n,k}/(n+1)]} \left| h(t) - h\left(\frac{I_{(i)}^{n,k}}{n+1}\right) \right|$$

$$+ \max_{0 \leq i \leq k} \left| \overline{h}_n\left(\frac{I_i^n}{n+1}\right) - h(U_i) \right| + \max_{1 \leq i \leq k} \left| h(U_{(i)}) - h\left(\frac{I_{(i)}^{n,k}}{n+1}\right) \right|,$$

where $U_{(i)}, 1 \leq i \leq k$, are the order statistics of $U_1, \ldots, U_k$. The desired convergence in probability is now a consequence of (31), (32), the fact that $I_k^n/(n+1)$ converges to $U_k$ a.s., and the a.s. continuity of $h$.

It was implicit in this proof that we were working with a strongly consistent family of discrete trees built from a self-similar fragmentation continuum tree and our usual argument shows that it still holds for any strongly consistent family. □

## 5. Beta-splitting, alpha and stable trees.

5.1. *Aldous's beta-splitting models.* Aldous [1] suggests a further study of what he calls *beta-splitting models*, where

$$\tilde{q}_n^{\text{Aldous}-\beta}(k) = \frac{1}{Z_n^{(\beta)}} \int_0^1 \binom{n}{k} x^{k+\beta} (1-x)^{n-k+\beta} \, dx$$

$$= \frac{1}{Z_n^{(\beta)}} \binom{n}{k} \frac{\Gamma(\beta+k+1)\Gamma(\beta+n-k+1)}{\Gamma(n+2\beta+2)},$$

$1 \leq k \leq n-1$, for some $-2 < \beta < \infty$. He says that these are sampling consistent and that he would like to establish continuum random tree limits (known only for $\beta = -3/2$, the Brownian CRT of Aldous [2]) also for all $-2 < \beta < -1$. He studies the asymptotic behavior of a randomly chosen leaf and heuristically argued that leaf-height functions rescaled in the same way should also converge. Our Theorem 2 and its height function ramification in Theorem 15 turn Aldous's heuristics into rigorous mathematics.

It is clear from Aldous's work [1] that the beta-splitting model with $-2 < \beta < -1$ corresponds to a binary dislocation measure

$$\nu_{\text{Aldous}-\beta}(s_1 \in dx) = C_\beta x^\beta (1-x)^\beta \mathbf{1}_{\{1/2 \leq x \leq 1\}}\, dx$$

and therefore satisfies the regular variation condition (6) with $\gamma = -\beta - 1$ and $\ell(x) \sim C_\beta/(-1-\beta)$. Since the splitting rules do not depend on $C_\beta$, we will choose $C_\beta = (-\beta - 1)/\Gamma(2+\beta)$ in the sequel.

Note that the symmetrized binary splitting rule above naturally gives rise to rooted *ordered* (or planar) trees $T_n^{\text{ord}}$ by the obvious recursive construction that builds tree $T_n^{\text{ord}}$ from a left subtree with $k$ leaves and a right subtree with $n-k$ leaves, with probability $\tilde{q}_n^{\text{Aldous}-\beta}(k)$, $1 \leq k \leq n-1$. We can now enumerate leaves from left to right and record their heights

$$h_n(i/(n+1)) = \text{distance from the root of the } i\text{th leaf from left to right}.$$
(33)

Also putting $h_n(0) = h_n(1) = 0$ and continuously extending to $[0,1]$ by linear interpolation gives the leaf-height function (which, in the binary case, fully encodes the discrete tree, just as the limiting height function fully encodes the limiting CRT) referred to by Aldous [1].

COROLLARY 16. *For a strongly sampling consistent family of trees $T_n^\circ$, $n \geq 1$, from the beta-splitting model with $-2 < \beta < -1$, we have*

$$\frac{T_n^\circ}{n^{-\beta-1}} \xrightarrow[n\to\infty]{(p)} \mathcal{T}_{(-\beta-1,\nu_{\text{Aldous}-\beta})}$$

*for the Gromov–Hausdorff metric. Furthermore, the associated rescaled leaf-height functions converge to an associated limiting height function* (*see Section* 4.3)

$$\left(\frac{h_n(t)}{n^{-\beta-1}}\right)_{0 \leq t \leq 1} \xrightarrow[n\to\infty]{(p)} (h_{-\beta-1,\nu_{\text{Aldous}-\beta}}(t))_{0 \leq t \leq 1}$$

*for the uniform norm.*

5.2. *Ford's alpha models.* There are several versions of the alpha model of random binary combinatorial trees, ordered and unordered, labeled and unlabeled, and each can be described in different ways; see Ford [21, 22]. We focus here on the induced distributions $P_n^{\text{ord},\circ}$ on $\mathbb{T}_n^{\text{ord},\circ}$, unlabeled shapes of planted (actually binary) *plane* (i.e., ordered) trees with $n$ leaves. Ford's original sequential construction leads to an increasing sequence of random trees $\tilde{T}_n^{\text{ord},\circ} \sim P_n^{\text{ord},\circ}$, $n \geq 1$, and we shall use this notation throughout our alpha model discussion. Fix $\alpha \in [0,1]$.

The sequential construction starts with the unique planted binary unlabeled plane trees $\tilde{T}_1^{\text{ord},\circ}$ and $\tilde{T}_2^{\text{ord},\circ}$ with one and two leaves, respectively. Given the random tree $\tilde{T}_n^{\text{ord},\circ}$ with $n$ leaves constructed following these rules, the $(n+1)$st leaf is added as follows: choose an edge according to weights $\alpha$ on edges between two inner vertices and $1-\alpha$ on edges between a leaf and an inner vertex. Since there are $n-1$ inner edges and $n$ leaf edges, the normalization constant is $n-\alpha$. Replace this edge between its two vertices by a new vertex and two edges linking its two vertices to the new vertex. Choose whether to attach the new leaf to the left or to the right of the new vertex with equal probability. The resulting random tree with $n+1$ leaves is called $\tilde{T}_{n+1}^{\text{ord},\circ}$.

We can now deduce the following corollary from Theorems 2 and 15.

COROLLARY 17. *Let $\tilde{T}_n^\circ$ be the unlabeled tree derived from Ford's sequential construction by forgetting the order of branches. Then*

$$\text{(34)} \qquad \frac{\tilde{T}_n^\circ}{n^\alpha} \xrightarrow[n\to\infty]{(d)} \mathcal{T}_{(\alpha,\nu_{\text{Ford}-\alpha})}$$

*for the Gromov–Hausdorff topology, where*

$$\nu_{\text{Ford}-\alpha}(s_1 \in dx)$$
$$= \frac{1}{\Gamma(1-\alpha)}(\alpha(x(1-x))^{-\alpha-1} + (2-4\alpha)(x(1-x))^{-\alpha})1_{\{1/2 \leq x \leq 1\}}\, dx.$$

*Furthermore, the associated rescaled leaf-height functions (33) encoding $\tilde{T}_n^{\text{ord},\circ}$ converge*

$$\left(\frac{h_n(t)}{n^\alpha}\right)_{0 \leq t \leq 1} \xrightarrow[n\to\infty]{(d)} (h_{\alpha,\nu_{\text{Ford}-\alpha}}(t))_{0 \leq t \leq 1}$$

*for the uniform topology on continuous functions defined on $[0,1]$.*

PROOF. Ford [21] shows that $(P_n^{\text{ord},\circ})_{n\geq 1}$ are the distributions of a sampling consistent Markov branching model with splitting kernel

$$\tilde{q}_n^{\text{Ford}-\alpha}(k) = \frac{\Gamma_\alpha(k)\Gamma_\alpha(n-k)}{\Gamma_\alpha(n)}\left(\frac{\alpha}{2}\binom{n}{k} + (1-2\alpha)\binom{n-2}{k-1}\right),$$
$$1 \leq k \leq n-1,$$

where $\Gamma_\alpha(n) = (n-1-\alpha)(n-2-\alpha)\cdots(2-\alpha)(1-\alpha) = \Gamma(n-\alpha)/\Gamma(1-\alpha)$.

For $0 < \alpha < 1$, Ford [22] also indicates that as $n \to \infty$, for all $0 < x < 1$,

$$n^{1+\alpha}\tilde{q}_n([xn]) \sim \frac{1}{\Gamma(1-\alpha)}\left(\frac{\alpha}{2}(x(1-x))^{-\alpha-1} + (1-2\alpha)(x(1-x))^{-\alpha}\right)$$
$$=: f_{\text{Ford}-\alpha}(x).$$

In the light of (16), we associate the binary dislocation measure

$$\nu_{\text{Ford}-\alpha}(s_1 \in dx) = (f_{\text{Ford}-\alpha}(x) + f_{\text{Ford}-\alpha}(1-x))1_{\{1/2 \le x \le 1\}}\,dx$$
$$= 2f_{\text{Ford}-\alpha}(x)1_{\{1/2 \le x \le 1\}}\,dx.$$

It is clear from Corollary 4 and the discussion which followed that the dislocation measure $\nu_{\text{Ford}-\alpha}$ induces Ford's splitting rule $(\tilde{q}_n)_{n\ge 2}$. By application of Theorem 2, (34) holds with $\xrightarrow{(d)}$ replaced by $\xrightarrow{(p)}$ for $T_n^\circ$ instead of $\tilde{T}_n^\circ$, where $(T_n^\circ)_{n\ge 1}$ is a *strongly sampling consistent* family derived from the homogeneous fragmentation with dislocation measure $\nu_{\text{Ford}-\alpha}$. But, according to Ford [21], for each fixed $n \ge 1$, there is the identity in distribution $\tilde{T}_n^\circ \sim T_n^\circ$. Theorem 15 can now be applied in the same way. □

As remarked by Ford [21], $\tilde{q}_\cdot^{\text{Ford}-\alpha} = \tilde{q}_\cdot^{\text{Aldous}-\beta}$ if and only if $\alpha = -\beta - 1 = 1/2$ (uniform model), $\alpha = \beta = 0$ (Yule model) or $\alpha = -\beta - 1 \uparrow 1$ (comb model). Also, we see that Ford's alpha model, as a model of exchangeable probability distributions on cladograms (by adding *exchangeable* leaf labels), is one of the wider class of Aldous's Markov branching models of type $c = 0$, $\tilde{\nu}(dx) = f(x)\,dx$ in Corollary 4.

Finally, we make some rather subtle points about Ford's sequential construction. It will be convenient to also consider $\tilde{T}_n^{\text{ord}}$ as the tree $\tilde{T}_n^{\text{ord},\circ}$ equipped with leaf labels in the order of Ford's sequential construction, and the unordered labeled tree $\tilde{T}_n$ derived from $\tilde{T}_n^{\text{ord}}$. In the following list, we consider $\alpha \in (0,1)$ and also exclude $\alpha = 1/2$, where no such subtleties arise.

- If a uniform leaf of $\tilde{T}_n^{\text{ord},\circ}$ is deleted, the tree generated by the remaining leaves has the same distribution as $\tilde{T}_{n-1}^{\text{ord},\circ}$. Nevertheless, for $\tilde{T}_n^{\text{ord}}$, with leaf labels in order of appearance, these labels are not exchangeable for $n \ge 3$. For example, in $\tilde{T}_3^{\text{ord}}$, leaf 3 has height 2 if the edge of $\tilde{T}_2^{\text{ord}}$ chosen for the insertion of 3 is adjacent to the root with probability $\alpha/(2-\alpha) \ne 1/3$.
- For fixed $n \ge 5$, the *joint* distribution of the unlabeled trees $(\tilde{T}_m^{\text{ord},\circ})_{1\le m\le n}$ is not the same as the joint distribution of $(\tilde{T}_n^{\text{ord},\circ,(m)})_{1\le m\le n}$, where $\tilde{T}_n^{\text{ord},\circ,(n)} = \tilde{T}_n^{\text{ord},\circ}$, and $\tilde{T}_n^{\text{ord},\circ,(m-1)}$ is obtained from $\tilde{T}_n^{\text{ord},\circ,(m)}$ by deleting a uniform leaf, $m = n, \ldots, 2$. Therefore, $(\tilde{T}_n^\circ)_{n\ge 1}$ is not strongly sampling consistent.
- We showed in Proposition 7 that for a strongly sampling consistent family of trees, convergence of finite-dimensional marginals holds almost surely with limiting trees $\mathcal{R}_k$ with edge lengths. In the next subsection, we will establish a corresponding result for $(\tilde{T}_n)_{n\geq 1}$. We also give a line-breaking construction of the almost sure limiting trees $\tilde{\mathcal{R}}_k$, $k \geq 1$.

- We conjecture that the completion of $\bigcup \tilde{\mathcal{R}}_k$ has the same distribution as $\mathcal{T}_{(\alpha,\nu_{\mathrm{Ford}-\alpha})}$ and that $\tilde{\mathcal{R}}_k$ can be embedded in $\mathcal{T}_{(\alpha,\nu_{\mathrm{Ford}-\alpha})}$ by suitable nonuniform, and presumably dependent, sampling of leaves.

5.3. *Limiting edge lengths in Ford's sequential construction.* Let $0 < \alpha < 1$. The limiting continuum random tree $\mathcal{T}_{(\alpha,\nu_{\mathrm{Ford}-\alpha})}$ naturally contains its uniformly sampled subtrees $\mathcal{R}_k$, $k \geq 1$, and ordered versions $\mathcal{R}_k^{\mathrm{ord}}$ are coded in $h_{\mathrm{Ford}-\alpha}$. If we denote the tree shape of $\mathcal{R}_n^{\mathrm{ord}}$ by $T_n^{\mathrm{ord}}$, then the distribution $P_n^{\mathrm{ord}}$ of $T_n^{\mathrm{ord}}$ is $P_n^{\mathrm{ord},\circ}$, equipped with *exchangeable* leaf labels. $\mathcal{R}_k^{\mathrm{ord}}$ is the almost sure scaling limits of the reduced trees $R(T_n^{\mathrm{ord}}, [k])$ as $n \to \infty$; see Proposition 7.

On the other hand, we naturally define $\tilde{P}_2^{\mathrm{ord}}$ to be uniform on the set $\mathbb{T}_2^{\mathrm{ord}}$ of two elements, and then $\tilde{P}_{n+1}^{\mathrm{ord}}$ directly from the sequential construction as the distribution of $\tilde{T}_{n+1}^{\mathrm{ord}}$, which is $\tilde{T}_n^{\mathrm{ord}}$ with the new leaf added according to Ford's rule and labeled $n+1$, that is, we label leaves in their order of appearance. In this setting, we also establish a.s. convergence of reduced subtrees.

PROPOSITION 18. (a) *For all $k \geq 1$, we have*

$$n^{-\alpha} R(\tilde{T}_n^{\mathrm{ord}}, [k]) \xrightarrow[n\to\infty]{\mathrm{a.s.}} \tilde{\mathcal{R}}_k^{\mathrm{ord}},$$

*in the sense of Gromov–Hausdorff convergence, where $(\tilde{\mathcal{R}}_k^{\mathrm{ord}})_{k\geq 1}$ is an increasing family of leaf-labeled $\mathbb{R}$-trees with edge lengths.*

(b) *The distribution of $\tilde{\mathcal{R}}_k^{\mathrm{ord}}$ is determined by the distributions of three independent random variables*: (i) *its shape* $\tilde{T}_k^{\mathrm{ord}} \sim \tilde{P}_k^{\mathrm{ord}}$; (ii) *its total length $S_k$ with density*

$$\frac{\Gamma(k+1-\alpha)}{\Gamma(k/\alpha)} s^{k/\alpha-1} g_\alpha(s),$$

*where $g_\alpha(s) = \frac{1}{\alpha} s^{-1-1/\alpha} f_\alpha(s^{-1/\alpha})$ is the Mittag–Leffler density derived from the stable density $f_\alpha$ with Laplace transform $e^{-\lambda^\alpha}$*; (iii) *Dirichlet edge length proportions $D_k = (D_k^{(1)}, \ldots, D_k^{(2k-1)}) \sim \mathcal{D}(1, \ldots, 1, (1-\alpha)/\alpha, \ldots, (1-\alpha)/\alpha)$, where, in $D_k$, we first list the $k-1$ inner edges, then the $k$ leaf edges, each by depth-first search.*

(c) $\tilde{\mathcal{R}}_k^{\mathrm{ord}}$ *is an inhomogeneous Markov process in its natural filtration $(\mathcal{H}_k)_{k\geq 1}$. More precisely, given $(\tilde{T}_k^{\mathrm{ord}}, S_k, D_k)$, the conditional distribution of $\tilde{T}_{k+1}^{\mathrm{ord}}$ is that where the Ford insertion happens at an edge $E_k$, sampled* *from the distribution on edges induced by $D_k$; $S_{k+1}$ has conditional density*

$$f_{S_{k+1}|S_k=z}(y) = \frac{\alpha^{1/\alpha}}{\Gamma((1-\alpha)/\alpha)}(y-z)^{1/\alpha-2}\frac{yg_\alpha(y)}{g_\alpha(z)};$$

*given $E_k$ is an inner edge, let $C_{k+1} \sim \mathrm{Unif}(0,1)$, otherwise, $C_{k+1} \sim \beta(1,(1-\alpha)/\alpha)$, independently from $S_{k+1}$; split $E_k$ into its proportions $C_{k+1}$ and $1-C_{k+1}$, $C_{k+1}$ being closer to the root. This determines the proportions $D_{k+1}$.*

PROOF. Fix $k \geq 1$ and $\tilde{T}_k^{\mathrm{ord}}$. For $n \geq k$, the reduced trees $R(\tilde{T}_n^{\mathrm{ord}},[k])$ all have the same shape as $\tilde{T}_k^{\mathrm{ord}}$. In the transition from $n$ to $n+1$, there may be no change of the reduced tree or one of the edge lengths may increase by 1. We can associate edges with $2k-1$ colors, where each edge in $\tilde{T}_k^{\mathrm{ord}}$ represents a color (but not white, which is reserved for later). Edges have weights which increase. Initially ($n=k$), the weights are one for each inner edge and $(1-\alpha)/\alpha$ for each leaf edge, zero for white. Each round, we pick a color at random, according to the current weights, and apply an updating rule as follows. Whenever an edge of the reduced tree is chosen (we recognize Ford's rule), the weight of that edge is increased by 1 and also the weight of white is increased by $(1-\alpha)/\alpha$. Whenever we pick white, the weight of white is increased by $1/\alpha$.

This model contains the essence of a Chinese restaurant process (see, e.g., [36], Lecture 3). Specifically, if we further discriminate the white weight by colored numbers identifying the subtree on the reduced tree in which the new leaf is added, then these subtrees can be considered tables in a restaurant and their leaves are customers. Suppose, at stage $n$, $m$ subtrees are present on $R(\tilde{T}_n^{\mathrm{ord}},[k])$. Each new customer joins any occupied table $i = 1, \ldots, m$ with probability $(n_i - \alpha)/(n - \alpha)$, where $n_i \geq 1$ is the number of customers already sitting at that table, and chooses a new table with remaining probability $(k + (m-1)\alpha)/(n - \alpha)$. This describes an $(\alpha, k - \alpha)$ seating plan in the terminology of [36].

(a)–(b) By [36], Theorem 3.8, the total number of tables scaled by $(n-k)^\alpha$ (where $n-k$ is the number of customers at stage $n$) converges almost surely so that for the total length $S_k^{(n)}$ of $R(\tilde{T}_n^{\mathrm{ord}},[k])$,

$$\frac{S_k^{(n)}}{n^\alpha} = \frac{S_k^{(n)} - 2k + 1}{(n-k)^\alpha} \frac{S_k^{(n)}}{S_k^{(n)} - 2k + 1} \frac{(n-k)^\alpha}{n^\alpha} \to S_k$$

and the distribution of $S_k$ is as specified.

In particular, if ignoring white, the total color weight still tends to infinity, even though it is asymptotically negligible against white weight. If we only record changes to the color weights, the restricted model still has the dynamics of the updating rule and so the pre-limiting proportions $D_k(n)$ converge

a.s. to the Dirichlet limit, as specified. Furthermore, $S_k^{(n)} D_k(n)/n^\alpha \to S_k D_k$ a.s. and, since the shape of reduced trees does not change (not even in the limit, as $D_k$ has only positive entries a.s.), this implies convergence in the Gromov–Hausdorff sense.

The independence of $\tilde{T}_k^{\mathrm{ord}}$, $S_k$ and $D_k$ can be seen by a conditioning argument: the independence of $(S_k, D_k)$ from $\tilde{T}_k^{\mathrm{ord}}$ follows since our argument actually gives us the conditional distribution of $(S_k, D_k)$ given $\tilde{T}_k^{\mathrm{ord}}$, which does not depend on $\tilde{T}_k^{\mathrm{ord}}$. Similarly, $S_k^{(n)}$ gives us the times at which the color weights change that leads to $D_k$; if we condition on $(S_k^{(n)})_{k \le n \le N}$, then we still observe the same dynamics of color weights and letting $N \to \infty$, we get independence of $D_k$ from the $\sigma$-field $\mathcal{S}_k$ generated by $(S_k^{(n)})_{n \ge k}$ with respect to which $S_k$ is measurable.

(c) Consider weight processes $\mathcal{W}_m$ leading to $\tilde{\mathcal{R}}_m^{\mathrm{ord}}$ as $1 \le m \le k+1$ varies. First, note that for $1 \le m < k$, $(\tilde{T}_m^{\mathrm{ord}}, S_m, D_m)$ is a measurable function of $(\tilde{T}_k^{\mathrm{ord}}, S_k, D_k)$. Therefore, the Markov property is trivially satisfied.

Now, let $\mathbf{t}_{k+1} \in \mathbb{T}_{k+1}^{\mathrm{ord}}$ be such that $k+1$ was added to an inner edge of the subtree $\mathbf{t}_k$ of $\mathbf{t}_{k+1}$, without loss of generality, directly to the left of the trunk. We then wish to calculate the expectation

$$\mathbb{E}(f(D_{k+1}, S_{k+1}) 1_{\{\tilde{T}_{k+1}^{\mathrm{ord}} = \mathbf{t}_{k+1}\}})$$

$$= \mathbb{P}(\tilde{T}_{k+1}^{\mathrm{ord}} = \mathbf{t}_{k+1}) \int \cdots \int f(e_1, \ldots, e_{2k}, 1 - e_1 - \cdots - e_{2k}, r)$$

$$\times \frac{\Gamma(k+1+k(1-\alpha)/\alpha)}{(\Gamma((1-\alpha)/\alpha))^k}$$

$$\times (e_{k+1} \cdots e_{2k}(1 - e_1 - \cdots - e_{2k}))^{(1-2\alpha)/\alpha}$$

$$\times \frac{\Gamma(k+2-\alpha)}{\Gamma((k+1)/\alpha)} r^{(k+1)/\alpha - 1} g_\alpha(r) \, d\mathbf{e} \, dr,$$

where $\mathbf{e} = (e_1, \ldots, e_{2k}, 1 - e_1 - \cdots - e_{2k})$ so as to identify the conditional distribution of $(\tilde{T}_{k+1}^{\mathrm{ord}}, C_{k+1}, S_{k+1})$ given $(\tilde{T}_k^{\mathrm{ord}}, D_k, S_k) = (\mathbf{t}_k, \mathbf{d}, s)$, where $\mathbf{d} = (d_1, \ldots, d_{2k-2}, 1 - d_1 - \cdots - d_{2k-2})$. We change variables

$$e_1 = \frac{d_1 c s}{r}, \qquad e_2 = \frac{d_1 (1-c) s}{r},$$

$$e_3 = \frac{d_2 s}{r}, \ldots, e_k = \frac{d_{k-1} s}{r}, \qquad e_{k+2} = \frac{d_k s}{r}, \ldots, e_{2k} = \frac{d_{2k-2} s}{r},$$

$$e_{k+1} = \frac{r - s}{r}$$

and calculate the Jacobian

$$\det\begin{pmatrix} \frac{cs}{r} & 0 & \cdots & 0 & \frac{d_1 s}{r} & \frac{d_1 c}{r} \\ \frac{(1-c)s}{r} & 0 & \cdots & 0 & -\frac{d_1 s}{r} & \frac{d_1(1-c)}{r} \\ 0 & \frac{s}{r} & \cdots & 0 & 0 & \frac{d_2}{r} \\ \vdots & \ddots & \ddots & \vdots & \vdots & \vdots \\ 0 & 0 & \cdots & \frac{s}{r} & 0 & \frac{d_{2k-2}}{r} \\ 0 & 0 & \cdots & 0 & 0 & -\frac{1}{r} \end{pmatrix} = \frac{d_1 s^{2k-1}}{r^{2k}}$$

by a development of the first row. This gives

$$\mathbb{E}(f(D_{k+1}, S_{k+1})1_{\{\tilde{T}^{\mathrm{ord}}_{k+1}=\mathbf{t}_{k+1}\}})$$
$$= \mathbb{P}(\tilde{T}^{\mathrm{ord}}_k = \mathbf{t}_k)$$
$$\times \int \cdots \int f\left(\frac{d_1 cs}{r}, \frac{d_1(1-c)s}{r}, \frac{d_2 s}{r}, \ldots, \frac{d_{k-1} s}{r}, \frac{r-s}{r}, \frac{d_k s}{r},\right.$$
$$\left.\ldots, \frac{d_{2k-2} s}{r}, \frac{(1-d_1-\cdots-d_{2k-2})s}{r}, r\right)$$
$$\times K(d_k \cdots d_{2k-2}(1-d_1-\cdots-d_{2k-2}))^{(1-2\alpha)/\alpha} s^{k/\alpha-1} g_\alpha(s)$$
$$\times d_1 s (r-s)^{(1-2\alpha)/\alpha} \frac{g_\alpha(r)}{g_\alpha(s)} d\mathbf{d}\, dc\, ds\, dr$$

for a positive constant $K$. We conclude that $\tilde{T}^{\mathrm{ord}}_{k+1}$, $C_{k+1}$ and $S_{k+1}$ are conditionally independent and that

$$\mathbb{P}(\tilde{T}^{\mathrm{ord}}_{k+1} = \mathbf{t}_{k+1} | \tilde{T}^{\mathrm{ord}}_k = \mathbf{t}_k, D_k = \mathbf{d}, S_k = s) = \tfrac{1}{2} d_1,$$
$$f_{C_{k+1}|\tilde{T}^{\mathrm{ord}}_k=\mathbf{t}_k, D_k=\mathbf{d}, S_k=s}(c) = 1,$$
$$f_{S_{k+1}|\tilde{T}^{\mathrm{ord}}_k=\mathbf{t}_k, D_k=\mathbf{d}, S_k=s}(r) = K_1 r(r-s)^{(1-2\alpha)/\alpha} \frac{g_\alpha(r)}{g_\alpha(s)}.$$

Similarly, if $k+1$ was added to a leaf edge of the subtree $\mathbf{t}_k$ of $\mathbf{t}_{k+1}$, without loss of generality, to the left of the first leaf edge in the order of depth-first search, which we may furthermore assume to be adjacent to the trunk in $\mathbf{t}_k$, then there will be an additional $(1-c)^{(1-2\alpha)/\alpha}$ in the change of variables since, now, $e_2$ and $e_{k+1}$ take the roles of $e_1$ and $e_2$, where $e_{k+1}$ is now the proportion of a leaf edge. We then get

$$\mathbb{P}(\tilde{T}_{k+1}^{\mathrm{ord}} = \mathbf{t}_{k+1} | \tilde{T}_{k}^{\mathrm{ord}} = \mathbf{t}_k, D_k = \mathbf{d}, S_k = s) = \tfrac{1}{2} d_k,$$

$$f_{C_{k+1}|\tilde{T}_{k}^{\mathrm{ord}}=\mathbf{t}_k, D_k=\mathbf{d}, S_k=s}(c) = K_2 (1-c)^{(1-2\alpha)/\alpha},$$

$$f_{S_{k+1}|\tilde{T}_{k}^{\mathrm{ord}}=\mathbf{t}_k, D_k=\mathbf{d}, S_k=s}(r) = K_1 r (r-s)^{(1-2\alpha)/\alpha} \frac{g_\alpha(r)}{g_\alpha(s)}. \quad \square$$

Proposition 18 is a generalization and refinement of [36], Exercises 7.4.10–7.4.13 dealing with the tree growth process in a Brownian excursion, $\alpha = 1/2$.

COROLLARY 19. *The counting process $N_t = \sup\{k \geq 0 : S_k \leq t\}$, $t \geq 0$, is a time-inhomogeneous renewal process with hazard function*

$$h_t(y) = \frac{y^{(1-\alpha)/\alpha - 1} t g_\alpha(t)}{\int_0^\infty (y+x)^{(1-\alpha)/\alpha - 1}(t+x) g_\alpha(t+x)\,dx},$$

*that is, the hazard rate is $h_t(y)$ at time $t$ if the last renewal occurred at $t - y \geq 0$.*

Note that since $f_{S_{m+1}|S_m = z}(y)$ integrates to 1, we have, for all $z \geq 0$,

$$\int_0^\infty y^{(1-\alpha)/\alpha - 1}(z+y) g_\alpha(z+y)\,dy = \frac{\Gamma((1-\alpha)/\alpha)}{\alpha} g_\alpha(z).$$

In the case $\alpha = 1/2$, we have $(1-\alpha)/\alpha - 1 = 0$ so that the tilting coefficients disappear and we can apply this formula to get $h_t(y) = t$. This is the Poisson line-breaking construction of the Brownian continuum random tree $\mathcal{T}_{\mathrm{Aldous}-(-3/2)} = \mathcal{T}_{\mathrm{Ford}-1/2}$ (see Aldous [2]), where the trees $\tilde{\mathcal{R}}_k \sim \mathcal{R}_k$ are constructed sequentially by breaking a line at the times of a Poisson point process in the wedge $\{(x,t) : t \geq 0, 0 \leq x \leq t\}$ with unit intensity per unit square. The heights of points generate the branch points on the previously grown tree.

Proposition 18(iii) can be interpreted as the line-breaking construction of the alpha model random tree. The inhomogeneous renewal process replaces the inhomogeneous Poisson arrival process at linearly increasing rate $t$. The branch points (heights of points in the point process) are no longer chosen uniformly as in the Brownian case, but with intensity skewed within each leaf edge, by the $\beta(1, (1-\alpha)/\alpha)$ choice replacing the uniform.

Denote by $V_1^{(n)}$ the number of leaves (out of $n$) added in (or as) subtrees to the left of the spine connecting leaf 1 to the root.

PROPOSITION 20. *We have a.s.*

$$\frac{V_1^{(n)}}{n} \to V_1 \sim \sum_{k=0}^\infty A_k W_k \left(\prod_{i=0}^{k-1} (1-W_i)\right),$$

where $W_i \sim \beta(1-\alpha, i\alpha+1-\alpha)$ are independent, $i \geq 0$, and, independently $A_k$, $k \geq 0$, are independent symmetric Bernoulli random variables.

PROOF. This is a consequence of the observation made in the proof of Proposition 18 that the partition of leaves according to subtrees is a Chinese restaurant. It is well known (see, e.g., [36]) that the table proportions are given by the products $W_m(1 - W_{m-1}) \cdots (1 - W_0)$. At the time of their creation, each subtree has an equal chance to grow on the left- and right-hand side of the spine, hence the result. $\square$

The distribution of $V_1$ is not a new distribution. It naturally arises in the more general context of size-biased sampling of Poisson point processes. Specifically, [35] identifies these atoms as the normalized jumps of a stable subordinator $\sigma$ with Laplace exponent $\lambda^\alpha$, tilted by $\sigma_1^{-(1-\alpha)}$, that is, we can also express the distribution of $V_1$ as

$$\mathbb{E}(f(V_1)) = \frac{\Gamma(2-\alpha)}{\Gamma(1/\alpha)} \mathbb{E}(f(\sigma_{1/2}/\sigma_1)\sigma_1^{-(1-\alpha)}).$$

See also [6, 36]. Recently, James, Lijoy and Pruenster [30] specified the density of $V_1$.

In general, $V_1$ does not have a uniform distribution as for $\alpha = 1/2$. Also, $V_1^{(k)}$ is not independent of $\tilde{T}_k^{\mathrm{ord}}$. For example, for $k = 3$, with two different shapes,

$$\mathbb{P}(V_1^{(3)} = 0 | \tilde{T}_3^{\mathrm{ord}} = {}^{\mathrm{v}}\mathrm{Y}) = \frac{1}{4-2\alpha} \quad \text{and} \quad \mathbb{P}(V_1^{(3)} = 0 | \tilde{T}_3^{\mathrm{ord}} = \mathrm{Y}^{\mathrm{v}}) = \frac{2-2\alpha}{4-2\alpha}$$

and these coincide if and only if $\alpha = 1/2$.

5.4. *Stable trees.* Duquesne and Le Gall [17] introduced a CRT that they called the stable tree $\mathcal{T}_{\mathrm{stable}-\alpha}^{\mathrm{ord}}$ of index $\alpha \in (1, 2]$, which describes the genealogy of a (continuous-state) stable branching process with a single infinitesimal ancestor conditioned to have unit total family size (integral of population sizes over time). For $\alpha = 2$, this is Aldous's Brownian continuum random tree, associated with Feller's diffusion. They have given the explicit distribution of the tree $R(\mathcal{T}_{\mathrm{stable}-\alpha}^{\mathrm{ord}}, L_1, \ldots, L_n)$ spanned by $n$ uniformly sampled leaves as follows. In fact, this identification of the finite-dimensional marginal distributions of $\mathcal{T}_{\mathrm{stable}-\alpha}^{\mathrm{ord}}$ may be taken as an alternative definition of the stable tree.

PROPOSITION 21 (Theorem 3.3.3 of [17]). (i) *Denote the shape of* $R(\mathcal{T}_{\mathrm{stable}-\alpha}^{\mathrm{ord}}, L_1, \ldots, L_n)$ *by* $T_n^{\mathrm{ord}}$. *Then*

$$\mathbb{P}(T_n^{\mathrm{ord}} = \mathbf{t}_n) = \frac{\alpha \Gamma(1-1/\alpha)}{\Gamma(n-1/\alpha)} \prod_{v \in \mathbf{t}_n, r_v \geq 2} \frac{(\alpha-1)\Gamma(r_v - \alpha)}{r_v! \Gamma(2-\alpha)}$$

for any $\mathbf{t}_n \in \mathbb{T}_n^{\mathrm{ord}}$, where $r_v$ is the number of children of vertex $v \in \mathbf{t}_n$.

(ii) *Given* $T_n^{\mathrm{ord}} = \mathbf{t}_n$, *the total length* $S_n$ *and the edge length proportions* $D_n$ *are conditionally independent;* $D_n$ *has a* $\mathcal{D}(1, \ldots, 1)$ *distribution on vectors of length* $l = |\mathbf{t}_n| - 1$, *the number of edges of* $\mathbf{t}_n$; $S_n$ *has density*

$$f_{S_n|T_n^{\mathrm{ord}}=\mathbf{t}_n}(s) = \frac{\alpha \Gamma(n-1/\alpha)}{\Gamma(\delta_{\mathbf{t}_n})\Gamma(l)}(\alpha s)^{l-1} \int_0^1 u^{\delta_{\mathbf{t}_n}-1} \eta(\alpha s, 1-u)\,du,$$

*where* $\delta_{\mathbf{t}_n} = n - 1/\alpha + (1 - 1/\alpha)l$ *and* $\eta(t, v)$ *is the density of a* $(1-1/\alpha)$-*stable subordinator* $(\sigma_t, t \geq 0)$ *with Laplace exponent* $\exp\{-\lambda^{1-1/\alpha}\}$.

There are a number of direct consequences.

COROLLARY 22 (Theorem 3.2.1 of [17], Lemma 5 of [33]). (i) *The tree shape without leaf labels,* $T_n^{\mathrm{ord},\circ}$, *is a Galton–Watson tree conditioned to have* $n$ *leaves, whose offspring distribution has probability generating function* $z + \alpha^{-1}(1-z)^\alpha$.

(ii) *The unordered tree shapes* $T_n^\circ$, $n \geq 1$, *form a strongly sampling consistent family of Markov branching models with splitting rule*

$$q_n(k_1, \ldots, k_r) = \frac{C_{k_1, \ldots, k_r} \Gamma(2-1/\alpha) \alpha^{-(r-2)} \Gamma(r-\alpha)}{\Gamma(n-1/\alpha)\Gamma(2-\alpha)} \prod_{j=1}^r \frac{\Gamma(k_j - 1/\alpha)}{\Gamma(1-1/\alpha)}$$

*for any* $r \geq 2$, $k_1 \geq \cdots \geq k_r \geq 1$, *where* $C_{k_1, \ldots, k_r}$ *is the combinatorial constant given in* (3).

Miermont [33, 34] studies fragmentation processes associated with $\mathcal{T}_{\mathrm{stable}-\alpha}$ and identifies the associated dislocation measure.

PROPOSITION 23 ([33]). *Let* $(\sigma_x, x \geq 0)$ *be a stable subordinator with Laplace exponent* $\lambda^{1/\alpha}$. *Denote by* $\Delta\sigma_{[0,1]} = (\Delta\sigma_x, x \in [0,1])^\downarrow$ *the jump sizes* $\Delta\sigma_x = \sigma_x - \sigma_{x-}$ *in decreasing order. Then* $\mathcal{T}_{\mathrm{stable}-\alpha}$ *is a* $(1-1/\alpha)$-*self-similar fragmentation CRT with dislocation measure*

$$\nu_{\mathrm{stable}-\alpha}(d\mathbf{s}) = \frac{\alpha^2 \Gamma(2-1/\alpha)}{\Gamma(2-\alpha)} \mathbb{E}\bigg(\sigma_1; \frac{\Delta\sigma_{[0,1]}}{\sigma_1} \in d\mathbf{s}\bigg).$$

*The associated Lévy measure* (18) *of the tagged particle subordinator is*

$$\Lambda_{\mathrm{stable}-\alpha}(dx) = \frac{\alpha-1}{\Gamma(1/\alpha)}(1-e^{-x})^{1/\alpha-2} e^{-(1-1/\alpha)x}\,dx.$$

By virtue of (20), which is equivalent to (6), the dislocation measure satisfies the regular variation condition with $\ell(x) \sim \alpha/\Gamma(1/\alpha)$ and also satisfies (7) for any $\rho > 0$ because the density of $\Lambda_{\mathrm{stable}-\alpha}$ decays exponentially as $x \to \infty$ (see also the discussion in Section 4.2). Therefore, we can apply Theorems 2 and 15.

COROLLARY 24. *For a strongly sampling consistent family of trees $T_n^\circ$, $n \geq 1$, from the Markov branching model with splitting rules identified in Corollary 22(ii) for some $1 < \alpha \leq 2$, we have*

$$\alpha \frac{T_n^\circ}{n^{1-1/\alpha}} \xrightarrow[n\to\infty]{(p)} \mathcal{T}_{\text{stable}-\alpha}$$

*for the Gromov–Hausdorff metric. Furthermore, the associated rescaled leaf-height functions (33) converge to the associated limiting height function (see Section 4.3)*

$$\left(\alpha \frac{h_n(t)}{n^{1-1/\alpha}}\right)_{0 \leq t \leq 1} \xrightarrow[n\to\infty]{(p)} (h_{1-1/\alpha, \nu_{\text{stable}-\alpha}}(t))_{0 \leq t \leq 1}$$

*for the uniform norm.*

It is known ([17]) that $\mathcal{T}_{\text{stable}-2} \sim 2\mathcal{T}_{\text{Aldous}-(-3/2)}$. Here, doubling a fragmentation CRT (i.e., all distances, or the associated height function) corresponds to halving the fragmentation rates. Also, for $\alpha \in (1,2)$, the factor $\alpha$ can be built into the limiting CRT as $\mathcal{T}_{\text{stable}-\alpha}/\alpha$, which is the CRT associated with the dislocation measure $\alpha\nu_{\text{stable}-\alpha}$.

Several papers study the convergence of conditioned discrete Galton–Watson trees. There are several different schemes of conditioning. The closest to our setting is conditioning on the total number of vertices in a Galton–Watson tree with offspring distribution in the domain of attraction of a stable law. Geiger and Kauffmann [23] study the convergence of reduced subtrees and show that the unconditional total length of $\mathcal{R}_k$ has a Gamma$(k + 1/\alpha, 1)$ distribution. Duquesne [16] establishes the convergence of associated height functions to the stable tree. It is not surprising that conditioning on the total number of leaves or the total number of vertices leads to the same limit, when suitably rescaled, since there are at most twice as many vertices as leaves and the ratio converges to $\alpha$ a.s.

Marchal [32] has a sequential construction of the shapes of the reduced stable tree similar to Ford's sequential construction of the alpha model. Marchal associates weights $1 - 1/\alpha$ with each edge, but also puts weight $k/\alpha - 1$ onto any vertex with $k$ subtrees. These weights also sum to $n - 1/\alpha$ at growth stage $n$. At each growth stage, an edge or vertex is chosen according to these weights and a new leaf edge added, either with an additional vertex in the "middle" of the edge or just attaching in the vertex increasing its number of subtrees by 1.

**Acknowledgments.** Thanks are due to Daniel Ford for discussing his alpha model with us at an early stage of his work. We would also like to thank two referees for valuable comments that led to an improvement of the presentation.

B. HAAS
CEREMADE
UNIVERSITÉ PARIS—DAUPHINE
PLACE DU MARÉCHAL DE LATTRE DE TASSIGNY
75775 PARIS CEDEX 16
FRANCE
E-MAIL: haas@ceremade.dauphine.fr

G. MIERMONT
DÉPARTEMENT DE MATHÉMATIQUES D'ORSAY
UNIVERSITÉ PARIS SUD 11
F-91405 ORSAY CEDEX
FRANCE
E-MAIL: miermont@math.u-psud.fr

J. PITMAN
DEPARTMENT OF STATISTICS
UNIVERSITY OF CALIFORNIA
BERKELEY, CALIFORNIA 94720
USA
E-MAIL: pitman@stat.berkeley.edu

M. WINKEL
DEPARTMENT OF STATISTICS
1 SOUTH PARKS ROAD
OXFORD OX1 3TG
UNITED KINGDOM
E-MAIL: winkel@stats.ox.ac.uk